\RequirePackage{fix-cm}
\documentclass{amsart} 
\usepackage[utf8]{inputenc}
\usepackage{mathptmx}
\usepackage{latexsym}
\usepackage{amsmath}
\usepackage{amssymb}
\usepackage{amsfonts} 
\usepackage{eucal} 
\usepackage{amsbsy}
\usepackage{amsthm}
\usepackage[all]{xy}
\usepackage{hyperref}

\newtheorem{thm}{Theorem}[section]
\newtheorem*{thm*}{Theorem}
\newtheorem{lemma}[thm]{Lemma}
\newtheorem{prop}[thm]{Proposition}
\newtheorem{cor}[thm]{Corollary}
\newtheorem*{cor*}{Corollary}

\theoremstyle{definition}
\newtheorem{defn}[thm]{Definition}
\newtheorem{example}[thm]{Example}

\theoremstyle{remark}
\newtheorem{remark}[thm]{Remark}

\newcommand {\hypco} {\ensuremath{\mbox{$\mathcal{H}$}}}
\newcommand {\dual}  {\ensuremath{\mathcal{D}}}
\newcommand {\dcp}   {\ensuremath{\mathbb{D}^\bullet}}

\newcommand {\La}    {\ensuremath{\mathcal{L}}}

\newcommand {\Pa}    {\ensuremath{\mbox{$\mathcal{P}$}}}

\newcommand {\Ua}    {\ensuremath{\mbox{$\mathcal{U}$}}}

\newcommand {\Wa}    {\ensuremath{\mathcal{W}}}
\newcommand {\Xa}    {\ensuremath{\mathcal{X}}}
\newcommand {\Ya}    {\ensuremath{\mathcal{Y}}}
\newcommand {\sha}   {\ensuremath{\textbf{A}^\bullet}}

\newcommand {\shq}   {\ensuremath{\textbf{Q}^\bullet}}

\newcommand {\shs}   {\ensuremath{\textbf{S}^\bullet}}

\newcommand {\sht}   {\ensuremath{\textbf{T}^\bullet}}

\newcommand {\der}   {\ensuremath{\mathbf{H}}}

\newcommand {\real}  {\ensuremath{\mathbb{R}}}
\newcommand {\intg}  {\ensuremath{\mathbb{Z}}}

\newcommand {\cplx}  {\ensuremath{\mathbb{C}}}
\newcommand {\rat}   {\ensuremath{\mathbb{Q}}}

\newcommand {\shom}  {\ensuremath{\operatorname{R \textbf{Hom}^\bullet}}}
\newcommand {\icm}   {\ensuremath{\textbf{IC}_{\bar{m}}^\bullet}}
\newcommand {\icn}   {\ensuremath{\textbf{IC}_{\bar{n}}^\bullet}}

\newcommand {\Lagr}  {\ensuremath{\operatorname{Lag}}}

\newcommand {\interi}{\operatorname{int}}

\newcommand {\omso}  {\ensuremath{\Omega_\ast^{\operatorname{SO}}}}

\newcommand {\smlhf} {\ensuremath{\mbox{$\frac{1}{2}$}}}
\newcommand {\zhf}   {\ensuremath{\intg [\smlhf]}}

\newcommand {\incl}  {\ensuremath{\operatorname{incl}}}
\newcommand {\cone}  {\ensuremath{\operatorname{cone}}}

\newcommand {\pt}    {\ensuremath{\operatorname{pt}}}
\newcommand {\id}    {\ensuremath{\operatorname{id}}}

\newcommand {\intr}   {\ensuremath{\operatorname{intr}}}
\newcommand {\osigma}   {\ensuremath{\overset{\circ}{\sigma}}}
\newcommand {\otau}   {\ensuremath{\overset{\circ}{\tau}}}

\newcommand {\OO}    {\ensuremath{\mathcal{O}}}
\newcommand {\Alb}    {\ensuremath{\operatorname{Alb}}}
\newcommand {\pr}  {\ensuremath{\mathbb{P}}}

\newcommand {\aff}    {\ensuremath{\mathbb{A}}}
\newcommand {\Ex}    {\ensuremath{\operatorname{Ex}}}
\newcommand {\Var}    {\ensuremath{\operatorname{Var}}}
\newcommand {\Spec}    {\ensuremath{\operatorname{Spec}}}

\begin{document}


\title[The L-Class of Covers and Varieties with Trivial Canonical Class]
    {Topological and Hodge L-Classes of Singular Covering Spaces and Varieties with Trivial Canonical Class}

\author{Markus Banagl}

\address{Mathematisches Institut, Universit\"at Heidelberg,
  Im Neuenheimer Feld 205, 69120 Heidelberg, Germany}

\email{banagl@mathi.uni-heidelberg.de}

\date{May, 2017}




\begin{abstract}
The signature of closed oriented manifolds is well-known to be 
multiplicative under finite covers. This fails for Poincar\'e complexes
as examples of C. T. C. Wall show. We establish the multiplicativity of
the signature, and more generally, the topological L-class, 
for closed oriented stratified pseudomanifolds that
can be equipped with a middle-perverse Verdier self-dual complex of sheaves,
determined by Lagrangian sheaves along strata of odd codimension
(so-called L-pseudomanifolds). This class of spaces contains all Witt spaces
and thus all pure-dimensional complex algebraic varieties.
We apply this result in proving the Brasselet-Schürmann-Yokura
conjecture for normal complex projective 3-folds with at most canonical singularities,
trivial canonical class and positive irregularity. The conjecture asserts the
equality of topological and Hodge L-class for compact complex algebraic rational
homology manifolds.
\end{abstract}

\maketitle


\tableofcontents


\section{Introduction}

For closed, oriented, smooth manifolds, the multiplicativity of the signature under finite covers
is a straightforward consequence of Hirzebruch's signature theorem, the naturality
of the cohomological L-classes of vector bundles and the multiplicativity of the
fundamental class under such coverings. For closed, oriented, topological manifolds,
this multiplicativity of the signature was proved by Schafer in \cite[Theorem 8]{schafer},
using Kirby-Siebenmann topological transversality.
On the other hand, for general Poincar\'e duality spaces which are not manifolds
multiplicativity fails. In fact, this phenomenon has been used by Wall to
construct Poincar\'e spaces that are not homotopy equivalent to manifolds.
In \cite{wall}, he constructs $4$-dimensional
Poincar\'e CW-complexes $X$ with cyclic fundamental
group of prime order such that the signature of both $X$ and its universal cover is $8$.
Such examples show that multiplicativity of the signature under finite covers
is a subtle matter in the presence of singularities.
In the present paper, we prove that the signature is multiplicative under
finite covers for the most general known class of pseudomanifolds allowing for
local Poincar\'e-Verdier self-duality and a bordism invariant signature.
These are the \emph{L-pseudomanifolds}, a class of pseudomanifolds that
includes all spaces that satisfy the Witt condition of Siegel \cite{siegel} and in particular
includes all equidimensional complex algebraic varieties.
Roughly, an oriented stratified pseudomanifold is an L-pseudomanifold if
it possesses a constructible Verdier self-dual complex of sheaves which interpolates
between lower and upper middle perversity intersection chain sheaves. (The latter
are dual to each other, but not self-dual unless the space satisfies the Witt condition.)
Every such complex of sheaves is given by Lagrangian subsheaves along the strata
of odd codimension. This class of spaces was introduced in \cite{banagl-mem} and is reviewed
here in Section \ref{sec.lpseudomanifolds}. 
In \cite{ablmp}, we related smoothly stratified L-pseudomanifolds and the self-dual perverse
sheaves on them to Cheeger's analytic Hodge $*$-invariant ideal boundary conditions,
\cite{cheeger1}, \cite{cheeger2}, \cite{cheeger3}.

Let $X$ be a closed L-pseudomanifold equipped with a Whitney stratification $\Xa$.
Homological characteristic L-classes $L_* (X;\Xa)\in H_* (X;\rat)$ were introduced in
\cite{banagl-lcl}, where we proved that the L-classes of
a choice of an interpolating self-dual sheaf complex as above are
independent of this choice, i.e. the choice of Lagrangian structures.
L-classes of general self-dual sheaves have first been considered
by Cappell and Shaneson in \cite{csstratmaps}; these do of course depend heavily on the choice
of self-dual sheaf. If $X$ is a Witt space, then $L_* (X;\Xa)$ agrees with the
Goresky-MacPherson-Siegel L-class of $X$ and if $X$ has only strata of
even codimension (in particular if $X$ is a complex algebraic variety), then
$L_* (X;\Xa)$ agrees with the Goresky-MacPherson L-class of \cite{gm1}, which,
following that reference, we shall often just denote as $L_* (X)$.
The first of two main results in this paper is the following multiplicativity statement
(Theorem \ref{thm.lclmult} in Section \ref{sec.multl}):

\begin{thm} (Multiplicativity of L-Classes for Finite Coverings.) \label{thm.lclassmultfincov}
Let $(X,\Xa)$ be a closed Whitney stratified L-pseudomanifold, 
$(X',\Xa')$ a Whitney stratified pseudomanifold
and $p: X' \to X$ a topological covering map of finite degree $d$ 
which is also a stratified map with respect to
$\Xa, \Xa'$. Then $(X', \Xa')$ is a closed L-pseudomanifold and
\[ p_* L_* (X',\Xa') = d\cdot L_* (X,\Xa), \]
where $p_*: H_* (X';\rat) \to H_* (X;\rat)$ is induced by $p$.
\end{thm}
Since for a space whose dimension is divisible by $4$,
the signature is the degree $0$-component of the L-class, it follows
that in particular the signature of L-pseudomanifolds (Witt spaces, equidimensional
complex varieties) is multiplicative under finite covers.
As far as our proof structure is concerned, we establish the signature statement
first and then use it as a basis for proving the L-class statement. 
In principle, there are ways to deduce such statements using ultimately the Atiyah-Singer
index theorem. Here, we wish to present proofs that do not require index theory.
Our chief tool is piecewise-linear signature homology $S^{PL}_* (-)$ as introduced in
\cite{banagl-ctslnw}. Based on M. Kreck's topological stratifolds as geometric cycles, 
and using topological transversality, topological signature
homology was first introduced by A. Minatta in \cite{minatta}. The coefficients had already
been introduced earlier in \cite[Chapter 4]{banagl-mem}.

In fact, Theorem \ref{thm.lclassmultfincov} can alternatively be deduced from
a more precise result on the L-class of a finite cover
(Theorem \ref{thm.vrr} in Section \ref{sec.multl}):
\begin{thm}
Let $(X,\Xa)$ be a closed $n$-dimensional
L-pseudomanifold, $(X',\Xa')$ an oriented Whitney stratified pseudomanifold
and $p: X' \to X$ an orientation preserving topological covering map of finite degree,
which is also a stratified map with respect to
$\Xa, \Xa'$. Then
\[ p_! L_* (X,\Xa) = L_* (X',\Xa'), \]
where $p_!: H_* (X;\rat) \to H_* (X';\rat)$ is the transfer induced by $p$.
\end{thm}

In the complex algebraic setting, results concerning the multiplicativity
of the $\chi_y$-genus (which corresponds to the signature for $y=1$)
under finite covers were obtained by A. Libgober and L. Maxim in
\cite[Lemma 2.3]{maximlibgober}. J. Schürmann discusses going up-and-down techniques for the
behavior of the motivic Chern class transformation $MHC_y$
under \'etale morphisms in \cite[Cor. 5.11, Cor. 5.12]{schurmann}.

In \cite{bsy}, Brasselet, Schürmann and Yokura conjectured that the Hodge-L-class $T_{1*} (X)$
and the topological Goresky-MacPherson L-class $L_* (X)$
are equal for compact complex algebraic varieties that are rational
homology manifolds.
We apply our Theorem \ref{thm.lclassmultfincov} in proving the 
Brasselet-Schürmann-Yokura conjecture for normal connected complex projective $3$-folds $X$
that have at worst canonical singularities, trivial canonical divisor, and  
$\dim H^1 (X;\OO_X)>0$.
The Hodge L-class $T_{1*} (X) \in H^{BM}_{2*} (X;\rat)$ in Borel-Moore homology is defined for
complex algebraic varieties (separated schemes of finite type over $\cplx$)
$X$ as the evaluation at $y=1$ of a Hirzebruch characteristic class
$T_{y*} (X) \in H^{\operatorname{BM}}_* (X)\otimes \rat [y]$.
That Hirzebruch class is the value on $[\id_X]$ of a Hirzebruch transformation
\[ T_{y*}: K_0 (\Var/X) \longrightarrow H^{BM}_{2*} (X)\otimes \rat [y], \]
where $K_0 (\Var /X)$ denotes the relative Grothendieck ring of complex algebraic
varieties over $X$ introduced by Looijenga in \cite{looijenga} in the context of motivic integration.
For nonsingular compact and equidimensional $X$, $T_{1*} (X) = L_* (X)$. 
Simple examples of curves show that generally $T_{1*} (X) \not= L_* (X)$ for singular $X$.
Let $\Omega (X)$ denote the abelian group of cobordism classes of
self-dual constructible complexes of sheaves on $X$.
For compact $X$, there is a commutative diagram of natural transformations
\[ \xymatrix{
K_0 (\Var/X) \ar[r]^{\operatorname{sd}} \ar[d]_{T_{y*}} & \Omega (X) \ar[d]^{L_*} \\
H_{2*} (X)\otimes \rat[y] \ar[r]^{y=1} & H_{2*} (X;\rat)
} \]
relating $T_{1*}$ to the aforementioned L-class transformation of Cappell-Shaneson.
But for general singular $X$, 
$\operatorname{sd} [\id_X]$ is not cobordant on $X$ to the intersection
chain sheaf of $X$.

Let $X$ be an equidimensional compact connected complex algebraic variety.
The degree zero component of $T_{y*} (X)$ is the Hodge polynomial $\chi_y (X)$.
As pointed out in \cite{cappmaxshan},
the natural transformation $MHT_y$ of \cite{bsy}, defined on the Grothendieck
group of algebraic mixed Hodge modules, can be evaluated on the intersection
chain sheaf, which, by Saito's work, carries a canonical mixed Hodge structure
\cite{saito89}, \cite{saito90}. This yields a homology characteristic class 
$IT_{y*} (X) = MHT_y [\mathbf{IC}^H (X)]$, whose associated genus is denoted by
$I\chi_y (X)$.
If $X$ is a rational homology manifold, then we have 
$\mathbf{IC}^H (X) \cong \rat^H_X \in D^b \operatorname{MHM}(X),$ and thus
$IT_{y*} (X)= MHT_y [\rat^H_X] =T_{y*} (X)$. In degree $0$ and at $y=1$, this means
$I\chi_1 (X) =\chi_1 (X)$. If $X$ is projective, then by Saito's 
intersection cohomology Hodge index theorem, $I\chi_1 (X)=\sigma (X)$.
Altogether then, $T_{1,0} (X)=\chi_1 (X)=\sigma (X)=L_0 (X)$, that is, in the
projective case the conjecture is known to hold in degree $0$. 
Furthermore, Cappell, Maxim, Schürmann and Shaneson \cite[Cor. 1.2]{cmssequivcharsing}
have shown that the conjecture holds for $X$ that are orbit spaces $Y/G$,
where $Y$ is a projective $G$-manifold and $G$ a finite group of algebraic
automorphisms.
The second main result of the present paper is:
\begin{thm} \label{thm.bsy}
The Brasselet-Schürmann-Yokura conjecture holds for normal, connected, projective, complex algebraic
$3$-folds $X$ with at worst canonical singularities, trivial canonical divisor
and $\dim H^1 (X;\OO_X)>0$: For such $X$,
\[ T_{1*} (X) = L_* (X). \]
\end{thm}
Our proof, which makes up Section \ref{sec.proofofbsy}, 
is predominantly geometric rather than Hodge-theoretic. 
We use the fact, reviewed in Section \ref{sec.canonicalsing}, that the Albanese map
is available for canonical singularities. When the canonical divisor is numerically
equivalent to zero, then a result of Kawamata provides splittings up to
a finite \'etale cover. That result is to some degree analogous to the
well-known decomposition theorem of Beauville-Bogomolov for compact Kähler manifolds with
trivial canonical class. In the case of canonical singularities, infinitesimal 
analogs of the Beauville-Bogomolov theorem, giving more information than Kawamata,
have been obtained by Greb, Kebekus and Peternell in \cite{grebetal}.
In addition to the multiplicativity of topological L-classes
under finite covers obtained in the present paper, 
we also use the fact that both the Hodge L-class and the
topological L-class of Witt spaces behave multiplicatively with respect to
cartesian products. 
On the Hodge-theoretic side, we apply a formula of
Cappell, Maxim, Schürmann and Shaneson,
\cite[Theorem 5.1]{cmssequivcharsing} (or alternatively
the Verdier Riemann-Roch type formula of
\cite[Cor. 3.1.3]{bsy}.)
Furthermore, we rely heavily on the additivity of the
Hodge L-class with respect to scissor relations and on ADE theory.
The proof shows that the $3$-folds under consideration are in fact 
rational homology manifolds. The most interesting case is 
$\dim H^1 (X;\OO_X)=1$, since it turns out that only in this case, nontrivial
abelian higher signatures and $L_2 (X)$ can occur. The proof
yields information on the range of those higher signatures.
We conclude by giving examples and L-class realization results based on constructions of
Reid \cite{reidfamous95}, Iano-Fletcher \cite{fletcher} and Alt{\i}nok \cite{altinokthesis}.\\

\textbf{Acknowledgements:} We thank Selma Alt{\i}nok-Bhupal,
Gavin Brown, Daniel Greb, Laurentiu Maxim and particularly Jörg Schürmann,
who suggested to bring in transfer maps, for helpful
communication and discussions.

\section{Pseudomanifolds Carrying a Bordism-Invariant Signature}
\label{sec.lpseudomanifolds}

Let $X^n$ be an $n$-dimensional oriented topological pseudomanifold with a fixed
topological stratification $\Xa$. The stratification specifies in particular a filtration
\[
X = X_n \supset X_{n-1} =X_{n-2} \supset X_{n-3} \supset \ldots \supset
X_0 \supset X_{-1} = \varnothing
\]
such that $X_j$ is closed in $X$
and every non-empty $X_j - X_{j-1}$ is an open topological manifold of dimension $j.$ 
A \emph{pure stratum} of $\Xa$ is a connected component of $X_j - X_{j-1}$.
We often write $\Sigma = X_{n-2}$ and call $X-\Sigma$ the \emph{top stratum} of $\Xa$.
Let $D(X,\Xa)$ denote the bounded derived category of complexes of
sheaves of real vector spaces on $X$, which are constructible with respect to $\Xa$.
Let $\dual$ denote the Verdier dualizing functor on $D(X,\Xa)$.
Regarding intersection sheaf complexes, we will follow the indexing conventions
of \cite{gm2}. Thus Verdier self-duality will be understood with respect to
$\dual [n]$, rather than with respect to $\dual$. 
Recall that $(X,\Xa)$ is said to satisfy the ($\real$-)\emph{Witt condition}, if
the canonical morphism $\icm (X) \to \icn (X)$ is an isomorphism, where
$\icm (X), \icn(X)$ denote the lower/upper middle perversity intersection chain sheaves
on $X$ with real coefficients. The condition holds trivially if there are no strata
of odd codimension. If the Witt condition is satisfied, then $\icm$ is a self-dual
sheaf, $\dual \icm (X) [n] \cong \icm (X)$ (in the derived category).
As shown in \cite{banagl-mem}, 
intersection homology type invariants on general topological pseudomanifolds 
that need not satisfy the Witt condition are given
by objects of a full subcategory $SD(X,\Xa)\subset D(X,\Xa),$ whose
objects satisfy two properties: On the
one hand, they are Verdier self-dual, on the other hand, they are
as close to the middle perversity intersection chain sheaves as
possible, that is, they interpolate between $\icm (X)$ and $\icn (X).$
We recall the precise definition:

\begin{defn} \label{def.sd}
Let $SD(X,\Xa)$ be the full subcategory of $D(X,\Xa)$ whose objects $\shs$
satisfy the following axioms: 
\begin{description}
 \item [(SD1)] Top stratum normalization: $\shs|_{X-\Sigma} \cong \real_{X-\Sigma} [n].$
 \item [(SD2)] Lower bound: $\der^i (\shs)=0,$ for $i<-n.$
 \item [(SD3)] Stalk condition for the upper middle perversity $\bar{n}:$ \\
  If $S$ is a pure stratum of $\Xa$ of codimension $k$, then 
  $\der^i (\shs |_S)=0$ for $i> \bar{n}(k)-n.$
 \item [(SD4)] Self-Duality: $\shs$ has an associated isomorphism 
 $d:\dual \shs [n] \stackrel{\cong}{\rightarrow} \shs$ \\
  such that
   $d$ is compatible with the orientation and $\dual d [n] = \pm d$.
\end{description}
\end{defn}
\noindent Here, $\der^i (\shs)$ denotes the cohomology sheaf of the 
complex $\shs$.
Depending on $(X,\Xa),$ the category $SD(X,\Xa)$ may or may not be empty.
If $\shs \in SD(X,\Xa),$ then there exist morphisms 
$\icm (X) \stackrel{\alpha}{\longrightarrow} \shs
\stackrel{\beta}{\longrightarrow} \icn (X)$ such that
\[ \xymatrix{
\icm (X)  \ar[r]^{\alpha} & \shs \\
\dual \icn (X)[n] \ar[u]_{\cong} \ar[r]^>>>>>{\dual \beta [n]} & \dual \shs [n] \ar[u]^{\cong}_d 
} \]
(where $d$ is given by $\bf{(SD4)}$) commutes, which clarifies the relation
between intersection chain sheaves and objects of $SD(X,\Xa).$
The main structure theorem on $SD(X,\Xa)$ is a description as a Postnikov
system with fibers given by
categories of Lagrangian structures along the strata
of odd codimension: 
\begin{thm} \label{thm.postn}
Let $n = \dim X$ be even. There is an equivalence of categories
\[ SD(X,\Xa) \simeq \Lagr (X_1 -X_0) \rtimes
  \Lagr (X_3 - X_2) \rtimes \ldots \rtimes
  \Lagr (X_{n-3} - X_{n-4}) \rtimes \operatorname{Coeff}(X-\Sigma). \]
(Similarly for $n$ odd.)
\end{thm}
While we will have no need to recall the definition of the Lagrangian categories
appearing on the right hand side of this equivalence, we will however illustrate
the theorem with examples.
\begin{example}
Let $X^6$ be the product of a circle with the (unreduced) suspension of
complex projective space, $X^6 = S^1 \times \Sigma \mathbb{C}P^2$,
equipped with the intrinsic stratification $\Xa$.
This stratification has a stratum of odd codimension $5$ consisting of the disjoint
union of two circles. The link of this stratum is $\mathbb{C}P^2$ and
there is no Lagrangian subspace in the middle cohomology 
$H^2 (\mathbb{C}P^2;\real)$ (e.g. the signature $\sigma (\mathbb{C}P^2) = 1\not= 0$).
Then the structure theorem implies $SD(X^6,\Xa)=\varnothing$, so that there
is no meaningful way to define locally self-dual intersection homology type invariants on
$X^6$. 

Let $X^4$ be the product of a circle with the suspension of
a torus, $X^4 = S^1 \times \Sigma T^2$, again equipped with the intrinsic stratification $\Xa$.
The stratum of odd codimension $3$ consists once more of 
two disjoint circles, but with link $T^2$.
There are many Lagrangian subspaces $\La$ in the middle cohomology 
$H^1 (T^2;\real)$, and
the structure theorem implies $SD(X^4,\Xa)\not= \varnothing$. In fact,
the functor from the right hand side of the equivalence of categories to $SD(X^4,\Xa)$
constructs a self-dual sheaf on $X^4$ for every choice of $\La.$
\end{example}

\begin{lemma} \label{lem.sdrefine}
If $\Xa$ and $\Xa'$ are topological stratifications of $X$ and $\Xa'$ is a refinement of $\Xa$,
then $SD(X,\Xa) \subset SD(X,\Xa')$ (full subcategory).
\end{lemma}
\begin{proof}
Let $\shs \in SD(X,\Xa)$. We have to check that $\shs$ is constructible
with respect to the stratification $\Xa'$ and
that it satisfies \textbf{(SD1)} -- \textbf{(SD4)} 
of Definition \ref{def.sd} with respect to $\Xa'$.

As for constructibility, let $S'$ be a pure stratum of $\Xa'$.
Since $\Xa'$ is a refinement of $\Xa$, there exists a pure stratum $S$
of $\Xa$ with $S' \subset S$. As $\shs$ is constructible with respect to $\Xa$,
the cohomology sheaf $\der^\bullet (\shs)|_S$ is locally constant.
The restriction of a locally constant sheaf to a subspace is locally constant.
Thus $\der^\bullet (\shs)|_{S'}$ is locally constant.

\textbf{(SD1)}: 
Let $S'$ be an $n$-dimensional pure stratum of $\Xa'$ and let
$S$ be a pure stratum of $\Xa$ such that $S' \subset S$.
Then $\dim S=n$ and thus $S$ is contained in the top stratum 
$X-\Sigma$ of $\Xa$. By \textbf{(SD1)} for $\Xa$,
$\shs|_S \cong \real_S [n]$ and thus
\[ \shs|_{S'} \cong (\shs|_S)|_{S'} \cong (\real_S [n])|_{S'} \cong \real_{S'} [n]. \]

\textbf{(SD2)}: The statement $\der^i (\shs)=0$ for $i<-n$ 
is independent of any choice of stratification.

\textbf{(SD3)}: 
Let $S'$ be a pure stratum of $\Xa'$ of codimension $k'$ and let
$S$ be a pure stratum of $\Xa$ such that $S' \subset S$.
Then the codimension $k$ of $S$ satisfies $k\leq k'$, which
implies $\bar{n}(k)\leq \bar{n}(k')$.
By \textbf{(SD3)} for $\Xa$, the cohomology sheaves
$\der^i (\shs|_S)$ vanish for $i>\bar{n}(k)-n$.
Thus 
$\der^i (\shs|_{S'}) = \der^i (\shs|_S)|_{S'} =0$
for $i>\bar{n}(k')-n$.

\textbf{(SD4)}: If $Y$ is any topological pseudomanifold, then the
dualizing complex $\dcp_Y$ is constructible with respect to any
topological stratification of $Y$. Thus $\dcp_X$ is constructible
with respect to $\Xa$ as well as $\Xa'$. Consequently, the dual
$\dual \shs = \shom (\shs, \dcp_X)$
is constructible with respect to both $\Xa$ and $\Xa'$, and
the isomorphism $\dual \shs [n] \cong \shs$ in $SD(X,\Xa)$ given by axiom 
\textbf{(SD4)} may be regarded as an isomorphism in $SD(X,\Xa')$. 
\end{proof}

Suppose that the dimension $n$ of $X$ is even, $n=2m$.
Given an object $\shs \in SD(X,\Xa),$ the duality isomorphism 
$d: \dual \shs [n] \cong \shs$ induces a nondegenerate
symmetric or anti-symmetric pairing
\[ \hypco^{-m} (X;\shs) \otimes \hypco^{-m} (X;\shs) \longrightarrow \real, \]
where $\hypco^* (X;\shs)$ denotes the global hypercohomology of $X$
with coefficients in $\shs$. 
We denote the signature of this pairing by $\sigma (X,\Xa;\shs)\in \intg$.
The first main result of \cite{banagl-lcl} is:
\begin{thm}(\cite[Thm. 4.1]{banagl-lcl}.)
Let $X^n$ be a closed oriented
topological pseudomanifold equipped with a topological stratification $\Xa$.
If $SD(X,\Xa)\not= \varnothing,$ then the signature
$\sigma (X,\Xa) = \sigma (X,\Xa;\shs) \in \intg$
is independent of the 
choice of $(\shs,d)$ in $SD(X,\Xa)$.
\end{thm}
By \cite[Thm. 4.1, p. 62]{banagl-mem},
this signature $\sigma (X,\Xa)=\sigma (X,\Xa;\shs)$ is a bordism invariant,
where admissible bordisms $(Y,\Ya)$ are compact oriented topologically stratified
pseudomanifolds-with-boundary,
whose interior $\interi Y$ can be covered by a complex of sheaves 
in $SD(\interi Y,\Ya|_{\interi Y})$ which restricts to the given boundary-sheaves near the boundary
(using collars).

For pseudomanifolds $X$ equipped with a Whitney stratification $\Xa$, homological
L-classes $L_i (X,\Xa;\sha)\in H_i (X;\rat)$
with coefficients in \emph{any} Verdier self-dual complex $\sha$ of sheaves
(not necessarily obeying the normalizations and costalk vanishing requirements
necessary for objects in $SD(X,\Xa)$) were defined by Cappell-Shaneson in
\cite{csstratmaps}. In general, these depend heavily on the choice of $\sha$, of course.
The second main result of \cite{banagl-lcl} is:
\begin{thm} \label{thm.independence}
Let $X^n$ be a closed, oriented
pseudomanifold equipped with a Whitney stratification $\Xa$.
If $SD(X,\Xa)\not= \varnothing,$ then the L-classes
$L_i (X,\Xa) = L_i (X,\Xa;\shs) \in H_i (X;\rat)$
are independent of the 
choice of $(\shs,d)$ in $SD(X,\Xa)$.
\end{thm}
\noindent Thus a Whitney stratified pseudomanifold has a well-defined 
homological L-class $L_* (X,\Xa),$ provided
$SD(X,\Xa)$ is not empty. If $X$ is a smooth closed oriented manifold
(equipped with any stratification $\Xa$), then $L_* (X,\Xa)$ agrees with the
Poincar\'e dual of the total Hirzebruch L-class of the tangent bundle of $X$.
If $(X,\Xa)$ satisfies the Witt condition, then $L_* (X,\Xa)$ agrees with the
Goresky-MacPherson-Siegel L-class.
It will be convenient to adopt the following terminology.
\begin{defn} \label{def.lpseudomfd}
An \emph{L-pseudomanifold}
is a Whitney stratified oriented pseudomanifold $(X,\Xa)$ such that
$SD(X,\Xa)$ is not empty.
\end{defn}
In our treatment of the BSY conjecture, we need the cartesian multiplicativity
of the topological L-class, at least in the Witt case. This has been
established by J. Woolf in \cite[Prop. 5.16]{woolf}.
(Note that the product of two Witt spaces is a Witt space.)
\begin{prop}(Woolf.) \label{prop.lclasscartesianmult}
Let $(X,\Xa),(Y,\Ya)$ be Whitney stratified oriented closed pseudomanifolds
both of which satisfy the Witt condition. Then the Goresky-MacPherson-Siegel
L-classes satisfy
\[ L_* (X\times Y,\Xa \times \Ya) = L_* (X,\Xa) \times L_* (Y,\Ya) \in H_* (X\times Y;\rat), \]
where $\Xa \times \Ya$ denotes the product stratification.
\end{prop}
Alternatively, this proposition can also be deduced from \cite{blm}, 
particularly Section 11 there.\\ 

Recall that a \emph{PL-space} is a topological space together with a class
of locally finite simplicial triangulations closed under linear subdivision and such
that any two triangulations in the class have a common subdivision.
All simplicial complexes in this paper will be locally finite.
An \emph{$n$-dimensional PL pseudomanifold}
is a PL-space $X$ such that for some (and hence every)
triangulation of $X$ by a simplicial complex $K$,
every simplex is contained in an $n$-simplex and
every $(n-1)$-simplex is a face of exactly two $n$-simplices.
An \emph{$n$-dimensional PL pseudomanifold-with-boundary} is a 
PL-space $Y$ such that for some (and hence every)
triangulation $L$ of $Y$
every simplex is contained in an $n$-simplex, and
every $(n-1)$-simplex is a face of either one or two $n$-simplices;
the union of the $(n-1)$-simplices which are faces of precisely one $n$-simplex is called
the \emph{boundary} of $Y$ and denoted $\partial Y$. One requires in addition that
the boundary $\partial Y$ is an $(n-1)$-dimensional pseudomanifold which 
is collared, that is, there is a PL embedding $\partial
Y\times [0,1)\to Y$ with open image which is the identity on $\partial Y$.
For a PL space $L$, let
$c^\circ L$ denote the open cone $([0,1)\times L)/(0\times L)$.  We recall the
inductive definition of a PL stratification of a PL pseudomanifold:
\begin{defn}
A \emph{PL stratification} $\Xa$ of a $0$-dimensional PL pseudomanifold $X$ is the filtration 
$X=X_0\supseteq X_{-1}=\varnothing$.
A \emph{PL stratification} $\Xa$ of an
$n$-dimensional PL pseudomanifold $X,$ $n>0,$
is a filtration by closed PL subspaces
\[ X=X_n \supset X_{n-1} = X_{n-2} \supset \cdots \supset X_0\supseteq
X_{-1}=\varnothing
\]
such that

\quad (a) $X-X_{n-1}$ is dense in $X$, and

\quad (b) for each point $x\in X_i - X_{i-1}$, there exists an open neighborhood
$U$ of $x$ in $X$ for which there is a  compact $(n-i-1)$-dimensional
PL stratified PL pseudomanifold  $(L,\Xa_L)$ and a PL homeomorphism
$\phi: \real^i \times c^\circ L\to U$
that takes $\real^i \times c^\circ (L_{j-1})$ onto $X_{i+j}\cap U$ for all $j$.
\end{defn}

The space $L$ in (b) is determined up to PL homeomorphism by $x$ and the
stratification $\Xa$; it is called the {\it link} of $X$ at $x$.
Since the cone on the empty set $L_{-1}$ is a point, taking $j=0$ in (b)
shows that $X_i - X_{i-1}$ is a manifold for every $i$.
A PL pseudomanifold $X$ always possesses a PL stratification $\Xa$.
Indeed, the \emph{intrinsic stratification} $\Xa_{\intr}$ is a PL stratification of $X$.
PL stratifications of PL pseudomanifolds-with-boundary are defined similarly;
the main point is that the collar PL homeomorphism is required to preserve 
strata.

Let $X$ be an $n$-dimensional PL pseudomanifold triangulated by a 
simplicial complex $K$.
The \emph{simplicial stratification $\Xa_K$ on $X$ induced by $K$} is given
by the filtration
\[ |K_n| = X \supset |K_{n-2}| \supset \cdots \supset |K_0|, \]
where $|K_i|$ is the underlying polyhedron of the $i$-skeleton $K_i$ of $K$, i.e. 
the union of all $i$-dimensional simplices of $K$.
The top stratum $X - |K_{n-2}|$ is the union of the interiors of
all $n$- and $(n-1)$-simplices of $K$. Then $\Xa_K$ is a PL stratification of $X$.
For every $K$, $\Xa_K$ is a refinement of the intrinsic stratification $\Xa_{\intr}$.
In particular, by Lemma \ref{lem.sdrefine},
\[ SD(X,\Xa_{\intr}) \subset SD(X,\Xa_K). \]

Every Whitney stratified pseudomanifold $(X,\Xa)$ can be triangulated 
by a simplicial complex $K$ so that
the strata of $\Xa$ are triangulated by subcomplexes, \cite{goresky}.
In this way, $X$ becomes a PL pseudomanifold such that the given
Whitney stratification $\Xa$ becomes a PL stratification of $X$ which is a
refinement of the intrinsic PL stratification $\Xa_{\intr}$. 
The simplicial stratification $\Xa_K$ is a refinement of $\Xa$ and thus
by Lemma  \ref{lem.sdrefine},
\[ SD(X,\Xa_{\intr}) \subset SD(X,\Xa) \subset SD(X,\Xa_K). \]
In particular, if $(X,\Xa)$ is an L-pseudomanifold, then $SD(X,\Xa_K)\not= \varnothing$.
In this situation, one may use simplicial maps to spheres and preimages
of points chosen in the interior of some top dimensional simplex of the sphere
to execute the Milnor-Pontrjagin-Thom construction. One obtains an L-class
$L_* (X,\Xa_K) \in H_* (X;\rat)$ which agrees with $L_* (X,\Xa)$
(\cite{csw},  \cite[\S 5]{csstratmaps}, \cite[\S 12.3]{weinberger}).

Let $Y$ be a topological space. Piecewise-linear signature homology $S_*^{PL} (Y)$ was
introduced in \cite{banagl-ctslnw}. Let us briefly recall its geometric cycle construction.
We describe first the objects that
represent elements in this theory. 
\begin{defn} \label{def.admissplrep}
An \emph{admissible PL-representative} is a quadruple
\[ (X, K, \shs, X\stackrel{f}{\longrightarrow} Y), \]
where
$X$ is an $n$-dimensional, closed, oriented PL pseudomanifold,
$K$ is a simplicial complex triangulating $X$,
$\shs \in SD(X, \Xa_K)$ is a self-dual complex of sheaves on $X$
(the self-duality isomorphism will be suppressed in the notation)
constructible with respect to the simplicial stratification $\Xa_K$, and
$f: X \rightarrow Y$ is a continuous map.
\end{defn}
Signature homology classes will then be bordism classes of admissible
PL-representatives as above. Thus we need to specify what a bordism between
such representatives is.
\begin{defn} \label{defn.admplbord}
A quadruple
\[ (W, L, \sht, W \stackrel{F}{\longrightarrow} Y), \]
is called an \emph{admissible PL-nullbordism} for an admissible PL-repre\-sentative
$(X^n, K, \shs, f)$ 
if:
\begin{itemize}
\item $W$ is a compact $(n+1)$-dimensional oriented PL pseudomanifold-with-boundary $\partial W=X$.
\item $L$ is a simplicial complex triangulating $W$.
\item $K\subset L$ is a subcomplex triangulating the subspace $X\subset W$.
\item $X$ is collared in $W$, i.e. there exists a closed PL neighborhood
 $N$ of $X$ in $W$ and an
 orientation preserving PL homeomorphism
 \[ c: X\times [0,1] \stackrel{\cong}{\longrightarrow} N \]
 such that $c(x,0)=x$ for all $x\in X$. We will use the notation
$\interi W = W-\partial W$.

\item $\sht \in SD(\interi W, \Wa_L)$, where $\Wa_L$ is the simplicial
 stratification induced on $\interi W$ by $L$.

\item Let $U$ be the open subset $U= c(X\times (0,1)) \subset W$ and let
 $c|$ denote the restriction $c|: X\times (0,1) 
 \stackrel{\cong}{\longrightarrow} U \subset \interi W$. The simplicial
 stratification on $W$ induces by restriction a stratification $\Ua$ of $U$.
 The simplicial stratification $\Xa_K$ on $X$ induces the product
 stratification $\Pa$ on $X\times (0,1)$. Then $c|$ is in general not
 stratum preserving with respect to $\Pa, \Ua$. We require:
 \begin{enumerate}
 \item $c|^\ast (\sht|_U)$ is $\Pa$-constructible on $X\times (0,1)$, and
 \item there is an isomorphism
  $c|^\ast (\sht|_U) \cong \pi^! \shs,$
  where $\pi: X\times (0,1) \rightarrow X$ is the projection.
 \end{enumerate}
\item $F: W \rightarrow Y$ is a continuous map such that $F|_X = f$. 
\end{itemize}
\end{defn}

Two admissible PL-representatives 
$(X_1, K_1, \shs_1, f_1)$ and $(X_2, K_2, \shs_2, f_2)$
are \emph{bordant} if there exists an admissible PL-nullbordism for
$(X_1, K_1, \shs_1, f_1) \sqcup -(X_2, K_2, \shs_2,$ $f_2)$.

\begin{defn}
\emph{Piecewise-linear signature homology} $S_n^{PL} (Y)$ is the set of
equivalence classes
\[ S_n^{PL} (Y) = \{ [(X^n,K,\shs, f)] ~|~ (S,K,\shs,f)
\text{ admissible PL-representative} \} \]
under the above bordism relation.
Disjoint union defines an abelian group structure on $S_n^{PL} (Y)$. 
\end{defn}
A continuous map $g: Y \rightarrow Y'$ induces a map
\[ S_n^{PL} (g): S_n^{PL} (Y) \longrightarrow S_n^{PL} (Y') \]
by
\[ S_n^{PL} (g)[(X,K,\shs,f)] = [(X,K,\shs, g \circ f)]. \]
Clearly $S_n^{PL} (1_Y) = 1_{S_n^{PL} (Y)}$ and
$S_n^{PL} (h\circ g) = S_n^{PL} (h) \circ S_n^{PL} (g)$.
Thus $S_\ast^{PL}$ is a functor from the category of topological spaces
and continuous maps to the category of abelian groups and group homomorphisms.
 In \cite{banagl-ctslnw}, we proved that
$S_\ast^{PL}$ is a homology theory on compact PL pairs.
The coefficients of piecewise linear signature homology are
\[ S_n^{PL} (\pt) \cong \begin{cases}
    \intg,~ n\equiv 0(4) \\
    0,~ n\not\equiv 0(4), \end{cases}
\]
where the isomorphism in the case $n\equiv 0(4)$ is given by the signature.
Let $\omso (-)$ denote bordism of smooth oriented manifolds.
A natural map
\[ 
\kappa_{PL} (Y): \omso (Y) \otimes_{\intg} \zhf \longrightarrow 
 S_\ast^{PL} (Y) \otimes_{\intg} \intg [\smlhf]
\]
is given on an element
$\text{$[$} M \stackrel{f}{\rightarrow} Y \text{$]$} 
\otimes a \in \omso (Y) \otimes_{\intg} \zhf$ as follows:
By J. H. C. Whitehead \cite{jhcwhiteh}, 
$M$ can be smoothly triangulated as a PL manifold by a simplicial complex $K$ and
this PL manifold is unique up to PL homeomorphism.
Define
\[ \kappa_{PL} (Y) (\text{$[$} M \stackrel{f}{\rightarrow} Y \text{$]$} 
\otimes a) = [(M, K, \real_{M} [n], 
 M \stackrel{f}{\rightarrow} Y)] \otimes a. \]
By \cite[p. 24]{banagl-ctslnw}, $\kappa_{PL} (Y)$ is a surjection.

\section{Multiplicativity of L-Classes Under Finite Coverings}
\label{sec.multl}

We prove that the signature and, more generally, the topological L-class
of an $L$-pseudo\-manifold (Definition \ref{def.lpseudomfd}) is multiplicative with respect
to covering projections of finite degree. Moreover, we demonstrate that
for such covers, the L-class of the covering space is the homology transfer of the
base's L-class, just as in the manifold situation.
We shall require the following standard concepts and facts from Whitney stratification theory.
Let $(X,\Xa)$ and $(Y,\Ya)$ be Whitney stratified spaces.
\begin{defn}
A proper continuous map $f:X\to Y$ is called a \emph{stratified map} with respect to
$\Xa, \Ya$ if there exist smooth manifolds $M,N$, embeddings
$X\subset M,$ $Y\subset N$ with Whitney's regularity conditions A and B satisfied
for the strata of $\Xa, \Ya$ (in particular, the pure strata are smooth submanifolds of $M,N$),
and an extension of $f$ to a smooth map $\widetilde{f}:M\to N$ such that
for every pure stratum $S$ in $\Ya$, $f^{-1} (S)$ is a union of connected
components of pure strata of $\Xa$, and
$f$ takes each of these components submersively to $S$.
\end{defn}
(Properness is required to obtain local topological triviality via Thom's first
isotopy lemma.)
This structure can always be achieved for proper complex algebraic maps:
\begin{prop} (\cite[p. 43]{gmsmt}) \label{prop.algwhitnstratmap}
Let $X,Y$ be complex algebraic subsets of complex algebraic smooth manifolds and
let $f:X\to Y$ be a proper complex algebraic map. Then there exist Whitney
stratifications of $X$ and $Y$ into pure strata that are complex algebraic
smooth manifolds,
such that $f$ becomes a stratified map.
\end{prop}

Subsequently, we shall often, without always mentioning this explicitly, make
use of the following simple but important observation, stating that if a covering projection is a
stratified map, then the stratification of the covering space must be pulled back from the base.
\begin{lemma} \label{lem.pullbackstrat}
Let $(X,\Xa),$ $(X',\Xa')$ be Whitney stratified spaces and $p:X' \to X$
a topological covering map which is a stratified map with respect to $\Xa, \Xa'$.
Then the pure strata in $\Xa'$ are the connected components of the preimages
$p^{-1} (S)$ of pure strata $S$ in $\Xa$.
\end{lemma}
\begin{proof}
Let $S'$ be a pure stratum in $\Xa'$. (Recall that we require pure strata to be connected.)
There exists a pure stratum $S$ in $\Xa$ such that $S' \cap p^{-1} (S) \not= \varnothing$.
This implies that $S' \subset p^{-1} (S)$, since $p$ is a stratified map.
Let $C \subset p^{-1} (S)$ be the connected component of $p^{-1} (S)$ that contains $S'$.
Suppose by contradiction that $S' \not= C$.
As $p$ is stratified, $p^{-1} (S)$ is a union of strata. Since they are all connected, $C$ is a 
union of strata. Thus $C$ is a (set-theoretically) disjoint union
\[ C = S' \cup \bigcup_{i\in I} S'_i,~ I\not= \varnothing,~ S'_i \not= \varnothing, \]
of strata of $\Xa'$.
The restriction $p|: p^{-1} (S)\to S$ is a covering map, and hence
$p|:C\to S$ is a covering map. It follows that $C$ is a topological manifold and
$\dim S' \leq \dim C = \dim S$. 
On the other hand, as $p$ is stratified, $p|: S' \to S$ is a smooth submersion so that
$\dim S' \geq \dim S$. We conclude that 
\[ \dim S' = \dim C = \dim S. \]
The same argument shows that for every $i\in I,$ $\dim S'_i = \dim C = \dim S'$.
This means that $S'$ and all $S'_i$ are codimension $0$ submanifolds (without boundary) of $C$,
and hence open subsets of $C$.
Therefore, $C-S' = \bigcup_{i\in I} S'_i$ is open, contradicting that $S'$ is connected.
We conclude that $S' =C$.
\end{proof}
We were unable to locate the following triangulation lemma in the literature and
thus present it here with full proof.
\begin{lemma} \label{lem.triangofcoverings}
Let $X$ be a topological space triangulated by a 
finite-dimensional
simplicial complex $K$ and let
$p:X' \to X$ be a topological covering map. Then there exists a triangulation of $X'$ by a simplicial
complex $K'$ such that $p$ is simplicial with respect to $K'$ and $K$.
\end{lemma}
\begin{proof}
Let $\tau: |K| \to X$ be the homeomorphism that provides the triangulation of $X$,
where $|K|$ denotes a geometric realization of the abstract simplicial complex $K$.
The triangulation $\tau$ can be lifted to a triangulation $\tau':|K'| \to X'$ of $X'$
such that for each simplex $\sigma'$ in $X'$,
$p|_{\sigma'}$ is a homeomorphism of $\sigma'$ onto a simplex $\sigma$ in $X$
(\cite[Thm. 24.6, p. 177]{moise}).
These homeomorphisms need not be linear, however. Thus $p$ will in general not
be simplicial with respect to $\tau, \tau'$. We shall now modify $\tau'$ (but not $\tau$)
inductively until $p$ becomes simplicial.

The map $p$ (using $\tau, \tau'$) induces a map $p_0: K'_0 \to K_0$ between the vertex sets.
We shall show first that $p_0$ gives rise to a morphism 
$p_*: K' \to K$ of abstract simplicial complexes. Thus we need to prove that
if vertices $v_0, \ldots, v_j$ in $K'$ span a simplex $\Delta'$ of $K'$, then the image vertices
$\tau^{-1} p\tau' (v_0), \ldots, \tau^{-1} p\tau' (v_j)$ span a simplex of $K$. 
Set $\sigma' = \tau' (\Delta')$ and $\sigma = p(\sigma')$.
Then $\sigma$ is a $j$-simplex in $X$ and there exist vertices
$w_0,\ldots, w_j$ in $K$ spanning a $j$-simplex $\Delta$ in $K$ with $\tau (\Delta)=\sigma$.
By induction on the dimension $j$, we shall prove the following claim:
There exists some permutation $\pi$ of the set $\{ 0,\ldots, j \}$ such that
$\tau (w_{\pi (i)}) = p\tau' (v_i)$ for all $i\in \{0,\ldots, j \}$.
For $0$-simplices $\Delta'$ this is clear by definition of the map $p_0$.
Suppose that $j>0$ and the claim holds for all $(j-1)$-simplices.
Let $\Delta'_1$ be the $(j-1)$-dimensional face of $\Delta'$ spanned by
$v_1,\ldots, v_j$ and let $\Delta'_{j-1}$ be the $(j-1)$-dimensional face of $\Delta'$ spanned by
$v_0,\ldots, v_{j-1}$. The images
$\sigma'_1 = \tau' (\Delta'_1),$ $\sigma'_{j-1} = \tau' (\Delta'_{j-1})$
lie in the boundary of $\sigma'$. As $p|:\sigma' \to \sigma$ is a homeomorphism,
it maps the boundary of $\sigma'$ onto the boundary of $\sigma$.
Therefore, $\sigma_1 = p(\sigma'_1)$ and $\sigma_{j-1} = p(\sigma'_{j-1})$
are $(j-1)$-simplices in $\partial \sigma$. So there exists an $i$ such that
$\sigma_{j-1} = \tau (\Delta_{j-1}),$ where $\Delta_{j-1}$ is the face of $\Delta$ in $K$
spanned by $w_0,\ldots, \widehat{w}_i,\ldots, w_j$.
Similarly, there exists a $k$ such that
$\sigma_1 = \tau (\Delta_1),$ where $\Delta_1$ is the face of $\Delta$ in $K$
spanned by $w_0,\ldots, \widehat{w}_k,\ldots, w_j$.
By induction hypothesis, there exists a bijection
\[ \beta_{j-1}: \{ 0,\ldots, j-1 \} \stackrel{\cong}{\longrightarrow}
 \{ 0, \ldots, \widehat{i},\ldots, j \} \]
such that $\tau (w_{\beta_{j-1} (l)}) = p\tau' (v_l)$ for all $l=0,\ldots, j-1$.
Similarly, there exists a bijection
\[ \beta_1: \{ 1,\ldots, j \} \stackrel{\cong}{\longrightarrow}
 \{ 0, \ldots, \widehat{k},\ldots, j \} \]
such that $\tau (w_{\beta_1 (l)}) = p\tau' (v_l)$ for all $l=1,\ldots, j$.
Define a map $\pi: \{ 0,\ldots, j \} \to \{ 0,\ldots, j \}$ by
\[ \pi (l) = \begin{cases}
\beta_{j-1} (l),& l\in \{ 0,\ldots, j-1 \}, \\
\beta_1 (l),& l=j.
\end{cases} \]
Let us verify that $\pi$ is a permutation:
We only need to see that $\beta_1 (j)=i$.
Suppose that $\beta_1 (j)\not= i$. Then $\beta_1 (j)\in \{ 0,\ldots, \widehat{i},\ldots, j \}$,
whence there is an $r\in \{ 0,\ldots, j-1 \}$ with
$\beta_{j-1} (r)=\beta_1 (j).$
It follows that
\[ p\tau' (v_r) = \tau (w_{\beta_{j-1} (r)}) =  \tau (w_{\beta_1 (j)}) = p\tau' (v_j). \]
Since $p|_{\sigma'}$ and $\tau'$ are homeomorphisms, $v_r = v_j$, i.e. $r=j$, a contradiction.
We conclude that $\beta_1 (j)=i$ as was to be shown.
Thus $\pi$ is indeed a permutation of $\{ 0,\ldots, j \}$ such that
$\tau (w_{\pi (l)}) = p\tau' (v_l)$ for all $l\in \{0,\ldots, j \}$ and the above claim is
established. Hence the vertex map $p_0$ induces a morphism $p_*: K' \to K$ of
abstract simplicial complexes.

Let $|p_*|: |K'| \to |K|$ be a geometric realization of $p_*$; this is a simplicial map.
We have shown above that $|p_*|$ has the following property:
Whenever $\sigma' = \tau' (\Delta')$ is a simplex in $X'$ with image simplex
$\sigma = p (\sigma') = \tau (\Delta)$ in $X$, then $|p_*| (\Delta') = \Delta$
and $|p_*|:\Delta' \to \Delta$ is a (simplicial) homeomorphism.

The construction of $\tau'$ ensures that $p$ is a cellular map, i.e.
$p (X'_j) \subset X_j,$ where $X'_j = \tau' (|K'|_j),$ $X_j = \tau (|K|_j),$
denote the $j$-skeleta. Thus we have restrictions
$p_j: X'_j \to X_j$ for all $j$. Note that $p_0: X'_0 \to X_0$ is already simplicial,
which furnishes the induction basis for an induction on the dimension $j$ of skeleta.

We shall now prove the following induction step for $j>0$:
Suppose that there is a triangulation $\tau':|K'| \to X'$ such that $p$ is
cellular with respect to $\tau', \tau$ and the restriction
$p_{j-1}: X'_{j-1} \to X_{j-1}$ is simplicial with respect to $\tau', \tau$, i.e.
the diagram
\[ \xymatrix{
|K'|_{j-1} \ar[r]^{\tau'|}_\cong \ar[d]_{|p_*|_{j-1}} & X'_{j-1} \ar[d]^{p_{j-1}} \\
|K|_{j-1} \ar[r]^{\tau|}_\cong & X_{j-1}
} \]
commutes, then
there exists a triangulation homeomorphism 
\[ \tau'':|K'| \to X' \] 
such that 
$\tau'' (|K'|_j) = X'_j =\tau' (|K'|_j)$ (and so $p$ is
cellular with respect to $\tau'', \tau$) and the restriction
$p_j: X'_j \to X_j$ is simplicial with respect to $\tau'', \tau$.

Let $\sigma'$ be a $j$-simplex of $K'$ in $X'$.
On $\sigma'$ we consider the map
$h (\sigma') = \tau \circ |p_*| \circ (\tau')^{-1}|$.
Writing $\sigma' = \tau' (\Delta')$ and
$\sigma = p (\sigma') = \tau (\Delta)$, where $\Delta, \Delta'$ are simplices
in $K,K'$ respectively, we have
$|p_*| (\Delta') = \Delta$ and $|p_*|:\Delta' \to \Delta$ is a (simplicial) homeomorphism
by the property of $|p_*|$ established above. It follows that
$h(\sigma')$ is a simplicial homeomorphism $h(\sigma'): \sigma' \to \sigma$.
Recall that the Alexander trick asserts that two homeomorphisms
$D^j \to D^j$ which agree on the boundary-sphere are isotopic rel boundary.
Thus the homeomorphism $p|:\sigma' \to \sigma$ is isotopic
rel $\partial \sigma'$ to the simplicial homeomorphism
$h (\sigma'): \sigma' \to \sigma.$
(Note that these two homeomorphisms indeed agree on $\partial \sigma'$ by
the commutativity of the above diagram.)
Let $\alpha_j: X'_j \to X'_j$ be the homeomorphism given on a $j$-simplex
$\sigma'$ by $(p|_{\sigma'})^{-1} \circ h(\sigma')$;
$\alpha_j$ is the identity on $X'_{j-1}$ and $\alpha_j (\sigma')=\sigma'$.
Then $p\circ \alpha_j: X'_j \to X_j$
is simplicial with respect to $\tau, \tau'$, since 
on simplices of dimension less than $j$ it is equal to $p_{j-1}$ (which is
simplicial by the above diagram) and on a $j$-simplex $\sigma'$, it is
given by
\[ (p\circ \alpha_j)|_{\sigma'} = p\circ (p|_{\sigma'})^{-1} \circ h(\sigma')
  = h(\sigma'), \]
which is simplicial by construction.

Now let $\sigma'$ be a $(j+1)$-dimensional simplex in $X'_{j+1}$.
The boundary $\partial \sigma'$ is a union of $j$-simplices and these are
simplexwise (but perhaps not pointwise) fixed by $\alpha_j$.
So $\alpha_j$ restricts to a homeomorphism
$\alpha_j|: \partial \sigma' \to \partial \sigma'$.
By the Alexander trick, this homeomorphism extends to a homeomorphism
$\sigma' \to \sigma'$. Doing this for every $\sigma'$, these extensions glue
to give a homeomorphism $\alpha_{j+1}: X'_{j+1} \to X'_{j+1}$ such that
$\alpha_{j+1} (\sigma')=\sigma'$ for all simplices $\sigma'$ in $X'_{j+1}$.
Continuing in this way we construct extending homeomorphisms 
$\alpha_{j+2}: X'_{j+2} \to X'_{j+2}$, etc., until we arrive at a
homeomorphism $\alpha = \alpha_n: X' \to X',$ 
where $n=\dim K' =\dim K < \infty$, which satisfies
\[ \alpha|_{X'_j} = \alpha_j,~ \alpha (X'_i) = X'_i \text{ for all } i. \]
Set 
\[ \tau'' = \alpha \circ \tau': |K'| \longrightarrow X'. \]
Then $\tau''$ is a homeomorphism such that
\[ \tau'' (|K'|_j) = \alpha \tau' (|K'|_j) = \alpha_j (X'_j) = X'_j \]
and the restriction
$p_j: X'_j \to X_j$ is simplicial with respect to $\tau'', \tau$, 
since on the $j$-skeleton of $|K'|,$ $\tau^{-1} \circ p \circ \tau''$ is given
by $\tau^{-1} \circ p\circ \alpha_j \circ \tau'$ and $p\circ \alpha_j$ is simplicial
with respect to $\tau, \tau'$. This finishes the proof of the induction step.
The desired triangulation of $X'$ is then $\tau''$ for $j$ with $X_j =X$.
\end{proof}

The next statement about simplicial maps is standard.
\begin{lemma} \label{lem.unionopensimplices}
Let $f:|K| \to |L|$ be a simplicial map between geometric realizations of simplicial
complexes $K,L$. Then the preimage $f^{-1} (\otau)\subset |K|$
of an open simplex $\otau$ in $|L|$ is a union of open simplices of $|K|$.
\end{lemma}

We are now in a position to discuss the behavior of the signature under 
finite coverings.
\begin{thm}  \label{thm.signmultfincov}
The signature of closed L-pseudomanifolds is multiplicative under finite coverings.
More precisely:
Let $(X,\Xa)$ be a closed L-pseudomanifold
whose dimension $n$ is a multiple of $4$, let
$(X',\Xa')$ be an oriented Whitney stratified pseudomanifold
and $p: X' \to X$ an orientation preserving topological covering map of finite degree $d$ 
which is also a stratified map with respect to
$\Xa, \Xa'$. Then $(X', \Xa')$ is a closed L-pseudomanifold and
\[ \sigma (X',\Xa')=d \cdot \sigma (X,\Xa). \]
\end{thm}
\begin{proof}
As $SD(X,\Xa)$ is not empty, we can and do choose a complex
$\shs \in SD(X,\Xa)$ with a Verdier self-duality isomorphism 
$d:\dual \shs [n] \cong \shs$.
Since the covering map $p$ is a local homeomorphism, the
sheaf functors $p^!$ and $p^*$ are canonically isomorphic.
Since moreover $p$ is a stratified map, $p^* \shs$ is constructible
with respect to $\Xa'$.
The isomorphism $d$ induces a self-duality isomorphism
\[ \dual (p^* \shs)[n] \cong p^! \dual \shs [n] \cong p^! \shs
  \cong p^* \shs \]
for the pullback $p^* \shs$. This proves axiom \textbf{(SD4)}
for $p^* \shs$.
For a general stratified map $f:X' \to X$ it need not be true that
$X' - \Sigma' \subset f^{-1} (X-\Sigma)$ for the top strata.
In the present case, however, $X' - \Sigma' \subset p^{-1} (X-\Sigma)$
holds, since $p$ is in addition a local homeomorphism.
As $\shs|_{X-\Sigma} \cong \real_{X-\Sigma} [n]$,
the pullback $p^* \shs$ is constant rank $1$ over $p^{-1} (X-\Sigma)$
and thus over $X' - \Sigma'$.
Hence \textbf{(SD1)} holds for $p^* \shs$.
The lower bound axiom clearly continues to hold for $p^* \shs$.
It remains to verify the
costalk vanishing axiom for $p^* \shs$.
Let $x' \in X'$ be any point and let $S'$ the unique pure stratum
containing $x'$. Let $k'$ be the codimension of $S'$.
Let $S$ be the unique pure stratum containing the image $p(x')$
and let $k$ denote the codimension of $S$.
Then, since $p$ is stratified, $S'$ is contained in $p^{-1} (S)$
and it follows that $k\leq k'$, using that $p$ is a local homeomorphism.
Thus $\der^i (p^* \shs)_{x'} = \der^i (\shs)_{p(x')}$ vanishes
for $i>k'$.
We have shown that $p^* \shs$ is in $SD(X',\Xa')$. 
In particular, this category is not empty and $(X',\Xa')$ is an
L-pseudomanifold.

Let $K$ be a simplicial complex triangulating $X$ such that
the Whitney strata of $\Xa$ are triangulated by subcomplexes (\cite{goresky}).
Then the simplicial stratification $\Xa_K$ of $X$ is a refinement 
of $\Xa$ and by Lemma \ref{lem.sdrefine},
$SD(X,\Xa) \subset SD(X,\Xa_K)$. We may thus regard
$\shs$ as an object of $SD(X,\Xa_K)$.

Let $K'$ be the triangulation of $X'$ provided by Lemma \ref{lem.triangofcoverings},
obtained by lifting $K$. Then $p:X' \to X$ is simplicial with respect to $K,K'$.
By Lemma \ref{lem.pullbackstrat}, the strata in $\Xa'$ are the connected components
of $p^{-1} (S)$ for strata $S$ in $\Xa$.
Thus the strata in $\Xa'$
are triangulated by subcomplexes of $K'$.
Let $\Xa_{K'}$ denote the simplicial stratification of $X'$ induced by $K'$.
Then $\Xa_{K'}$ refines $\Xa'$, $SD(X',\Xa') \subset SD(X',\Xa_{K'})$, and
\[ \sigma (X',\Xa') = \sigma (X',\Xa'; p^* \shs) = \sigma (X',\Xa_{K'}; p^* \shs). \]
The second equality holds because hypercohomology does not see 
stratifications, it only sees the sheaf $p^* \shs$.
The quadruple $(X,K,\shs, \id_X)$ is an admissible PL-representative
in the sense of Definition \ref{def.admissplrep} and hence defines a class
\[  [(X,K,\shs, \id_X)] \in S_n^{PL} (X).     \]
As the canonical map
\[ \kappa_{PL} (X): \omso (X) \otimes_{\intg} \zhf 
 \longrightarrow 
 S_\ast^{PL} (X) \otimes_{\intg} \zhf \]
is surjective,
there exists a closed, smooth, oriented $n$-manifold $M$ and a continuous
map $f: M\to X$ such that for a smooth triangulation $L$ of $M$,
\[ [(M,L, \real_M [n], f)] = 2^s [(X,K,\shs, \id_X)] \in S_n^{PL} (X),  \]
for some nonnegative integer $s$.
Consequently, there is an admissible PL-bordism
\[ (W^{n+1}, Q, \sht, W \stackrel{F}{\longrightarrow} X), \]
where
$W$ is a compact oriented $(n+1)$-dimensional PL pseudomanifold-with-boundary
triangulated by the simplicial complex $Q$,
$\partial W = 2^s X \sqcup -M,$ $F|_{\partial W} = 2^s \id_X \sqcup f,$
$2^s K \sqcup L$ is a subcomplex of $Q$ triangulating $\partial W$ and
$\sht$ is a self-dual sheaf in $SD(\interi W, \Wa_Q)$, where $\Wa_Q$ is the
simplicial stratification of the interior of $W$ induced by $Q$.
Furthermore, we know that there is a PL collar
$c:\partial W \times [0,1] \to N$ onto a closed polyhedral neighborhood $N$ of
$\partial W$ in $W$ such that
\begin{equation} \label{equ.collart}
c|^* (\sht|_U) \cong \pi^! (2^s \shs \sqcup \real_M [n]), 
\end{equation}
where $c|$ is the restriction of $c$ to $\partial W \times (0,1)$, $U$ is
the image of this open cylinder under $c$, and
$\pi: \partial W \times (0,1)\to \partial W$ is the first factor projection.
Let $p_W: W' = W\times_X X' \to W$ be the pullback covering of degree $d$,
fitting into the cartesian diagram
\[ \xymatrix{
W' \ar[r]^{F'} \ar[d]_{p_W} & X' \ar[d]^p \\
W \ar[r]^F & X.
} \]
Using Lemma \ref{lem.triangofcoverings}, triangulate $W'$
by a simplicial complex $Q'$ such that $p_W$ is simplicial with respect
to $Q,Q'$. The space $W'$ is a compact oriented $(n+1)$-dimensional PL
pseudomanifold-with-boundary, 
\[
\partial W' 
= p_W^{-1} (\partial W) = p_W^{-1} (2^s X \sqcup -M) 
= p_W^{-1} (2^s X) \sqcup -p_W^{-1} (M).
\]
The first preimage, $p_W^{-1} (2^s X),$ is in fact $2^s X',$
and the restriction of $F'$ to $2^s X'$ is $2^s \id_{X'}$.
Since $K'$ has been obtained by lifting $K,$ and $2^s K$ is a subcomplex of $Q$,
the complex $2^s K'$ is a subcomplex of $Q'$.
Restricting the above square to $M' = p_W^{-1} (M),$ we obtain a pullback
diagram of degree $d$ covering spaces
\begin{equation} \label{equ.cartdiammprxxpr}
 \xymatrix{
M' \ar[r]^{f'} \ar[d]_{p_M} & X' \ar[d]^p \\
M \ar[r]^f & X.
} 
\end{equation}
The boundary of $W'$ can then be written as $\partial W' = 2^s X' \sqcup -M'$.
Restricting the triangulation of $W'$ to $M'$, we obtain a triangulation of $M'$
by a simplicial subcomplex $L'$ of $Q'$. Then $p_M$ is simplicial with respect to
$L,L'$, as it is the restriction of $p_W$, which is simplicial with respect to $Q,Q'$.
By \cite[Prop. 2.12]{lee}, $M'$ has a unique smooth structure such
that $p_M$ is a smooth covering map and a local diffeomorphism. 
Let $p_0: W' - \partial W' \to W-\partial W$
be the restriction of $p_W$ to the interior and set
$\shq = p_0^* \sht$.
We claim that
\[  (W', Q', \shq, W' \stackrel{F'}{\longrightarrow} X') \]
is an admissible PL-nullbordism for
\[ 2^s (X', K', p^* \shs, \id_{X'}) \sqcup 
  -(M', L', \real_{M'} [n], M' \stackrel{f'}{\longrightarrow} X'). \]
It remains to be shown that 
$\shq \in SD(\interi W', \Wa_{Q'})$, where $\Wa_{Q'}$ is the simplicial
stratification of the interior of $W'$ induced by $Q'$, and that
$\shq$ looks like $2^s p^* \shs \sqcup \real_{M'} [n]$ near the boundary.
Let $c:\partial W \times [0,1] \cong N$ be the collar used in (\ref{equ.collart}).
The complex $\sht$ is constructible with respect to $\Wa_Q$ and hence
$\der^i (\sht)|_{\osigma}$ is constant for every open simplex $\osigma$ of $Q$.
Thus $\der^i (p_0^* \sht)$ is constant on $p^{-1}_0 (\osigma)$ for every $\sigma \in Q$
which does not lie in the boundary of $W$.
Since $p_W$ is simplicial,
Lemma \ref{lem.unionopensimplices} implies that $p^{-1}_W (\osigma)$ is a 
union of open simplices $\otau$ of $Q'$.
It follows that if $\otau$ is the interior of any simplex $\tau \in Q'$ not in $\partial W'$,
then $\der^i (p_0^* \sht)|_{\otau}$ is constant.
We conclude that the pullback $p_0^* \sht$ is constructible over the $(n-2)$-skeleton of
$\Wa_{Q'}$. But it is also constructible over the top stratum of $\Wa_{Q'}$, using
that this top stratum is contained in the preimage under $p_0$ of the top stratum
of $\Wa_Q$, as $p_W$ is simplicial and a local homeomorphism.
The arguments put forth at the beginning of this proof then
show that $p_0^* \sht$ is an object of $SD(\interi W', \Wa_{Q'})$.
The space $N' = p_W^{-1} (N)$ is a polyhedral neighborhood of 
$\partial W'$ in $W'$.
Let $\partial W_0$ be a connected component of $\partial W$ and
let $\partial W'_0$ be a connected component of $p^{-1}_{\partial W} (\partial W_0)$.
Then the restriction of $p$ defines a covering projection
$p_0: \partial W'_0 \to \partial W_0$.
Let $N_0$ be the connected component of $N$ such that $\partial W_0 \subset N_0$
and let $N'_0$ be the connected component of $N'$ such that $\partial W'_0 \subset N'_0$.
Then the restriction of $p$ defines a covering projection
$p_N: N'_0 \to N_0$. Another covering of $N_0$ is given by the composition
\[ \partial W'_0 \times [0,1] \stackrel{p_0 \times \id}{\longrightarrow}
 \partial W_0 \times [0,1] \stackrel{c|}{\longrightarrow} N_0. \]
Since $c$ is the identity on $\partial W_0 \times \{ 0 \},$ the groups
$\pi_1 (\partial W'_0) = \pi_1 (\partial W'_0 \times [0,1])$ and
$\pi_1 (N'_0)$ have conjugate images in $\pi_1 (N_0)$.
Therefore, there exists a homeomorphism $c'_0: \partial W'_0 \times [0,1] \to N'_0$
such that
\[ \xymatrix{
\partial W'_0 \times [0,1] \ar[r]_{c'_0}^\cong \ar[d]_{p_0 \times \id} & N'_0 \ar[d]^{p_N} \\
\partial W_0 \times [0,1] \ar[r]_c^\cong & N_0
} \]
commutes. Doing this for all connected components, we arrive at a collar
$c':\partial W' \times I \cong N'$ so that
\[ \xymatrix{
\partial W' \times [0,1] \ar[r]_{c'}^\cong \ar[d]_{p_{\partial W} \times \id} & N' \ar[d]^{p_W|} \\
\partial W \times [0,1] \ar[r]_c^\cong & N
} \]
commutes, where $p_{\partial W}$ denotes the restriction of $p_W$ to
the boundary.
With $U' = c' (\partial W' \times (0,1)) = p_W^{-1} (U),$ this diagram restricts to
\[ \xymatrix{
\partial W' \times (0,1) \ar[r]_{c'|}^\cong \ar[d]_{p_{\partial W} \times \id} & U' \ar[d]^{p_W|} \\
\partial W \times (0,1) \ar[r]_{c|}^\cong & U.
} \]
Using the standard factor projection $\pi': \partial W' \times (0,1) \to \partial W',$ the diagram
\[ \xymatrix{
\partial W' \times (0,1) \ar[r]^{\pi'} \ar[d]_{p_{\partial W}\times \id} & \partial W' \ar[d]^{p_{\partial W}} \\
\partial W \times (0,1) \ar[r]^{\pi}  & \partial W
} \]
commutes. Thus, using the isomorphism (\ref{equ.collart}),
\begin{align*}
c'|^* (\shq|_{U'})
&= c'|^* (p_0^* \sht|_{U'}) = c'|^* (p_W|)^* (\sht|_U) \\
&= (p_{\partial W} \times \id)^* c|^* (\sht|_U) \cong
     (p_{\partial W} \times \id)^* \pi^! (2^s \shs \sqcup \real_M [n]) \\
&= (p_{\partial W} \times \id)^! \pi^! (2^s \shs \sqcup \real_M [n])
 = \pi'^! p_{\partial W}^! (2^s \shs \sqcup \real_M [n]) \\
&= \pi'^! p_{\partial W}^* (2^s \shs \sqcup \real_M [n])
  = \pi'^! (2^s p^* \shs \sqcup p_M^* \real_M [n]) \\
&= \pi'^! (2^s p^* \shs \sqcup \real_{M'} [n]),
\end{align*}
which establishes the claim.
Consequently,
\[ 2^s [(X', K', p^* \shs, \id_{X'})] = 
  [(M', L', \real_{M'} [n], f')] \in S_n^{PL} (X'). \]
Under the composition
\[ S_n^{PL} (X') \longrightarrow S_n^{PL} (\pt) 
 \stackrel{\sigma}{\longrightarrow} \intg, \]
this element maps to the signature
$2^s \sigma (X',\Xa') = \sigma (M').$
By multiplicativity of the signature for the smooth covering $p_M: M' \to M,$
\[ \sigma (M') = d\cdot \sigma (M). \]
Under the composition
\[ S_n^{PL} (X) \longrightarrow S_n^{PL} (\pt) 
 \stackrel{\sigma}{\longrightarrow} \intg, \]
the element
\[ 2^s [(X,K,\shs, \id_X)] = [(M,L, \real_M [n], f)] \in S_n^{PL} (X)  \]
maps to $2^s \sigma (X,\Xa) = \sigma (M)$
and thus
\[ \sigma (X',\Xa') = 2^{-s} \sigma (M') = 2^{-s} d\sigma (M) = d\sigma (X,\Xa). \]
\end{proof}

\begin{cor} \label{cor.signmultwitt}
The signature of closed Whitney stratified Witt spaces is multiplicative under finite coverings.
More precisely:
Let $(X,\Xa)$ be a closed oriented Whitney stratified pseudomanifold
whose dimension $n$ is a multiple of $4$ and which satisfies the Witt condition.
Let
$(X',\Xa')$ be an oriented Whitney stratified pseudomanifold
and $p: X' \to X$ an orientation preserving topological covering map of finite degree $d$ 
which is also a stratified map with respect to
$\Xa, \Xa'$. Then $(X', \Xa')$ satisfies the Witt condition and
the Goresky-MacPherson-Siegel signatures satisfy
\[ \sigma (X')=d \cdot \sigma (X). \]
\end{cor}
\begin{proof}
The Witt condition ensures that $\icm (X) \cong \icn (X)$ is a self-dual complex
of sheaves such that $\icm (X) \in SD(X,\Xa)$. Then 
\[ \icm (X') \cong p^* \icm (X) \cong p^* \icn (X) \cong \icn (X') \]
is in $SD(X',\Xa')$ and $(X',\Xa')$ satisfies the Witt condition.
Note also that the GMS signatures
$\sigma (X) = \sigma (X,\Xa,\icm (X)),$
$\sigma (X') = \sigma (X',\Xa',\icm (X'))$ are stratification independent 
by topological invariance of intersection homology.
\end{proof}

\begin{cor} \label{cor.signmultalg}
The Goresky-MacPherson signature of equidimensional, possibly singular, normal
complex projective varieties is multiplicative under algebraic finite \'etale covers.
\end{cor}
\begin{proof}
Let $p: X' \to X$ be an algebraic finite \'etale cover of irreducible 
projective normal complex algebraic
varieties. By Proposition \ref{prop.algwhitnstratmap}, 
$X'$ and $X$ can be equipped with Whitney stratifications
$\Xa'$ and $\Xa$ such that $p$ becomes a stratified map.
The strata are complex algebraic manifolds, so in particular have even dimension
and the Witt condition for $(X,\Xa)$ is vacuously satisfied.
Furthermore, $p$ is indeed a topological covering map (see e.g. 
\cite[Rem. 3.5]{grebetaletalefund}).
Equidimensionality ensures that $X$ and $X'$ are pseudomanifolds. 
The statement thus follows from Corollary \ref{cor.signmultwitt}.
\end{proof}

The above signature multiplicativity statements will now constitute the basis
of more general multiplicativity statements for L-classes.

\begin{thm} (Multiplicativity of L-Classes for Finite Coverings.) \label{thm.lclmult}
Let $(X,\Xa)$ be a closed $n$-dimensional
L-pseudomanifold, $(X',\Xa')$ an oriented Whitney stratified pseudomanifold
and $p: X' \to X$ an orientation preserving topological covering map of finite degree $d$ 
which is also a stratified map with respect to
$\Xa, \Xa'$. Then $(X', \Xa')$ is a closed L-pseudomanifold and
\[ p_* L_* (X',\Xa') = d\cdot L_* (X,\Xa), \]
where $p_*: H_* (X';\rat) \to H_* (X;\rat)$ is induced by $p$.
\end{thm}
\begin{proof}
Let $j: Y^m \hookrightarrow X^n$ be a normally nonsingular inclusion
of an oriented, closed, stratified pseudomanifold $Y^m.$ Consider an open
neighborhood $E\subset X$ of $Y,$ the total space of an $\real^{n-m}$-
vector bundle over $Y,$ and put $E_0 = E-Y,$ the total space with the
zero section removed. Let $u\in H^{n-m} (E,E_0;\rat)$ denote the Thom class.
With $\pi: E\rightarrow Y$ the projection, the composition
\[ H_k (X;\rat) \rightarrow H_k (X,X-Y;\rat)
  \overset{e_\ast}{\underset{\cong}{\leftarrow}} H_k (E,E_0;\rat)
  \overset{u\cap -}{\underset{\cong}{\rightarrow}} H_{k-n+m} (E;\rat)
  \overset{\pi_\ast}{\underset{\cong}{\rightarrow}} H_{k-n+m} (Y;\rat) \]
defines a Gysin map
\[ j^!: H_k (X;\rat) \longrightarrow H_{k-n+m} (Y;\rat). \]
Let $\epsilon_*: H_0 (Y;\rat)\to \rat$ be the augmentation map.
By \cite[Prop. 8.2.11]{banagl-tiss}, the L-class $L_* (X,\Xa;\sht)$ of any self-dual
$\Xa$-constructible bounded complex $\sht$ of sheaves is uniquely determined by
\[ \epsilon_* j^! L_{n-m} (X, \Xa; \sht) = \sigma (Y, \Ya; j^! \sht) \] 
for every normally nonsingular inclusion $j: Y^m \hookrightarrow
X^n$ with trivial normal bundle, where the stratification $\Ya$ is obtained by intersecting
$Y$ with the strata in $\Xa$.
Let $\shs$ be an object in $SD(X,\Xa)$. Then $p^* \shs$ is an object of $SD(X',\Xa')$
and the pushforward $Rp_* p^* \shs$ remains self-dual, since $p$ is a proper map.
Let $d\shs$ denote the direct sum of $d$ copies of $\shs$. We shall prove that
\begin{equation} \label{equ.ldspps}
L_* (X,\Xa; d\shs) = L_* (X,\Xa; Rp_* p^* \shs) \in H_* (X;\rat). 
\end{equation}
By the above uniqueness statement, it suffices to show that
\begin{equation} \label{equ.signjdsjpps} 
\sigma (Y,\Ya; j^! d\shs) = \sigma (Y,\Ya; j^! Rp_* p^* \shs) 
\end{equation}
for every normally nonsingular inclusion $j: Y^m \hookrightarrow
X^n$ with trivial normal bundle.
Now 
\[ \sigma (Y,\Ya; j^! d\shs) = d\sigma (Y,\Ya; j^! \shs) = d\sigma (Y), \]
as $j^! \shs \in SD(Y,\Ya)$.
Let $Y' = p^{-1} (Y)$ and consider the cartesian diagram
\[ \xymatrix{
Y' \ar[r]^{j'} \ar[d]_{p_Y} & X' \ar[d]^p \\
Y \ar[r]^j & X
} \]
Then $p_Y$ is a topological covering map of degree $d$.
We shall show that $p_Y$ is a stratified map with respect to 
$\Ya, \Ya',$ where the stratification $\Ya'$ of $Y'$ is obtained by intersecting with the
strata of $\Xa'$.
Let $S$ be a pure stratum in $\Ya$. We need to verify that $p^{-1}_Y (S)$ is a union
of strata of $\Ya'$. Let $S'$ be a stratum of $\Ya'$ with $S' \cap p^{-1}_Y (S)\not= \varnothing$.
By definition of $\Ya, \Ya'$ there exist pure strata $T$ in $\Xa$ and $T'$ in $\Xa'$ such that
$S=T\cap Y,$ $S' = T' \cap Y'$. As $S' \cap p^{-1}_Y (S)$ is contained in
$T' \cap p^{-1} (T),$ the latter intersection is nonempty. Thus $T' \subset p^{-1} (T)$,
as $p$ is stratified with respect to $\Xa, \Xa'$.
So if $y' \in S'=T'\cap Y'$, then $p(y')\in T\cap Y=S$, which shows that
$S' \subset p^{-1}_Y (S)$.
It remains to see that $p_Y|: S' \to S$ is a submersion.
The restriction $p|: p^{-1} (T)\to T$ is a topological covering map.
By Lemma \ref{lem.pullbackstrat}, $T'$ is a connected component of $p^{-1} (T)$.
Hence $p_T = p|:T' \to T$ is still a covering map. 
Furthermore, $p_T$ is a smooth submersion, since $p$ is stratified with respect to 
$\Xa, \Xa'$. By the inverse function theorem, $p_T$ is a local diffeomorphism.
Let $y' \in S'$ be any point and $y=p_Y (y')\in S$.
Suppose that $\gamma: (-\epsilon, \epsilon)\to S$ is a smooth curve in $S$ with
$\gamma (0)=y$. Let $U\subset T$ be an open neighborhood of $y$ such that there is
a smooth inverse $p^{-1}_T: U \stackrel{\cong}{\longrightarrow} U',$
with $U' \subset T'$ an open neighborhood of $y'$. We may assume that $\epsilon >0$
is so small that $\gamma (-\epsilon, \epsilon) \subset U\cap S$.
Let $\gamma': (-\epsilon, \epsilon) \to U'$ be the smooth curve 
$\gamma' = p^{-1}_T \circ \gamma$. Then $\gamma' (0)=y'$ and if $t\in (-\epsilon, \epsilon),$
then $\gamma' (t) \in Y' = p^{-1} (Y)$, since
\[ p\gamma' (t) = p p^{-1}_T \gamma (t) = \gamma (t)\in S \subset Y. \]
Thus the image of $\gamma'$ lies in $S'$ and we have $p_Y \gamma' = \gamma$.
This shows that the derivative of $p_Y|:S' \to S$ at $y'$ is onto.
We have shown that $p_Y$ is indeed a stratified map with respect to $\Ya, \Ya'$.

According to \cite[Prop. V.10.7(4)]{borel},
\[ j^! Rp_* (p^* \shs) = Rp_{Y*} j'^! (p^* \shs). \]
Hence,
\[ \sigma (Y,\Ya; j^! Rp_* p^* \shs) = \sigma (Y', \Ya'; j'^! p^* \shs) = \sigma (Y'), \]
as $j'^! p^* \shs \in SD(Y',\Ya')$.
By the multiplicativity of the signature for $L$-pseudomanifold coverings,
Theorem \ref{thm.signmultfincov},
\[ \sigma (Y')=d \sigma (Y), \]
establishing (\ref{equ.signjdsjpps}), and thus also (\ref{equ.ldspps}).
By \cite[Prop. 5.4]{csstratmaps}, 
\[ L_* (X,\Xa; d\shs) = dL_* (X,\Xa; \shs) = dL_* (X;\Xa), \]
while by \cite[Prop. 5.5]{csstratmaps},
\[ L_* (X,\Xa; Rp_* p^* \shs) = p_* L_* (X', \Xa'; p^* \shs) = p_* L_* (X', \Xa'). \]
\end{proof}

\begin{cor} \label{cor.lclassmultwittcover}
The Goresky-MacPherson-Siegel L-class of Whitney stratified Witt spaces 
is multiplicative under finite coverings.
More precisely:
Let $(X,\Xa)$ be an oriented closed Whitney stratified pseudomanifold
which satisfies the Witt condition.
Let
$(X',\Xa')$ be an oriented Whitney stratified pseudomanifold
and $p: X' \to X$ an orientation preserving topological covering map of finite degree $d$ 
which is also a stratified map with respect to
$\Xa, \Xa'$. Then $(X', \Xa')$ satisfies the Witt condition and
the Goresky-MacPherson-Siegel L-class satisfies
\[ p_* L_* (X',\Xa')=d \cdot L_* (X,\Xa) \in H_* (X;\rat). \]
\end{cor}

\begin{cor}  \label{cor.lclassmultalgcover}
The Goresky-MacPherson L-class of equidimensional, possibly singular, normal
complex projective varieties is multiplicative under algebraic finite \'etale covers.
\end{cor}

Our arguments above can be rephrased, and strengthened, in terms of transfers 
to yield a Verdier-Riemann-Roch type theorem. 
Given a covering $p: X' \to X$ of finite degree $d$, there is a transfer 
\[ p_!: H_* (X;\rat) \to H_* (X';\rat) \]
such that $p_* p_!$ is multiplication by $d$. Furthermore,
given a cohomology class $\xi \in H^* (X;\rat)$ and a homology class
$x\in H_* (X;\rat)$, the cap formula
\begin{equation} \label{equ.transfernat}
p_! (\xi \cap x) = (p^* \xi) \cap p_! (x)
\end{equation}
holds. Given a cartesian diagram
\[ \xymatrix{
Y' \ar[r]^{f'} \ar[d]_q & X' \ar[d]^p \\
Y \ar[r]^f & X,
} \]
there is a base change relation
\begin{equation} \label{equ.basech}
f'_* q_! = p_! f_*.
\end{equation}
(For further information on transfers, see for example \cite[Ch. V, Sect. 6]{boardman}.)
\begin{lemma} \label{lem.vrrsmooth}
If $p: M' \to M$ is a finite smooth covering of smooth, oriented, closed manifolds, then
\[ p_! L_* (M) = L_* (M'). \]
\end{lemma}
\begin{proof}
The transfer of the fundamental class is $p_! [M] = [M']$.
Let $TM, TM'$ be the tangent bundles of $M,M'$.
Using the identity (\ref{equ.transfernat}), we calculate:
\begin{align*}
p_! L_* (M)
&= p_! (L^* (TM)\cap [M]) = (p^* L^* (TM))\cap p_! [M] \\
&= L^* (p^* TM) \cap [M'] = L^* (TM')\cap [M'] = L_* (M').
\end{align*}
\end{proof}

\begin{lemma} \label{lem.flbordinvonsignhom}
Let $Y$ be a topological space. Then the assignment
\[ (X^n, K, \shs, f) \mapsto f_* L_* (X,\Xa_K) \]
induces a well-defined natural transformation
\[ S^{PL}_n (Y) \longrightarrow H_* (Y;\rat). \]
\end{lemma}
\begin{proof}
Let $(W^{n+1},L,\sht,F:W\to Y)$ be an admissible PL-nullbordism for
$(X^n, K,\shs, f:X\to Y)$. 
We must show that $f_* L_* (X,\Xa_K)=0$.
The construction of \cite[Section 2]{bcs} and the proof of Proposition 1 in \emph{loc. cit.}
carries over from Whitney stratified Witt spaces to Whitney or PL stratified pseudomanifolds
$(X,\Xa)$ with $SD(X,\Xa)\not= \varnothing$, using $\shs \in SD (X,\Xa)$ in place of
the intersection chain sheaf $\icm (X)$.
Thus the nullbordism $W$ has a relative L-class
$L_* (W,\Wa_L) \in H_* (W,X;\rat),$ where $\Wa_L$ denotes the simplicial stratification
of $W$ induced by the complex $L$, such that
$\partial_* L_* (W,\Wa_L) = L_{*-1} (X,\Xa_K)$ under the connecting homomorphism
in the exact sequence of the pair $(W,X)$,
\[ H_* (W,X) \stackrel{\partial_*}{\longrightarrow} H_{*-1} (X)
 \stackrel{j_*}{\longrightarrow} H_{*-1} (W), \]
where $j: X=\partial W \hookrightarrow W$ is the inclusion of the boundary.
Since $f= Fj,$ we have
\[ f_* L_* (X,\Xa_K) = F_* j_* L_* (X,\Xa_K) = F_* j_* \partial_* L_{*+1} (W,\Wa_L)=0. \] 
To prove naturality, let $g:Y\to Y'$ be a continuous map.
Then the diagram
\[ \xymatrix{
S^{PL}_n (Y) \ar[r] \ar[d]_{S^{PL}_n (g)} & H_* (Y;\rat) \ar[d]^{g_*} \\
S^{PL}_n (Y') \ar[r]  & H_* (Y';\rat)
} \]
commutes, as both the clockwise and the counterclockwise composition send
$[(X^n,K,\shs,f)]$ to $g_* f_* L_* (X,\Xa_K)\in H_* (Y';\rat)$.
\end{proof}
For the special case of Witt spaces, the above proof yields a natural
L-class transformation $\Omega^{\operatorname{Witt}}_n (Y)\to H_* (Y;\rat)$
on Witt bordism. This transformation has already been described by
Woolf \cite[Cor. 5.15, Prop. 5.16]{woolf} and by 
Curran \cite[Lemma 4.8]{curran}.
We can now give a precise computation of the L-class of a finite covering of
a singular space:
\begin{thm} \label{thm.vrr}
Let $(X,\Xa)$ be a closed $n$-dimensional
L-pseudomanifold, $(X',\Xa')$ an oriented Whitney stratified pseudomanifold
and $p: X' \to X$ an orientation preserving topological covering map of finite degree,
which is also a stratified map with respect to
$\Xa, \Xa'$. Then
\[ L_* (X',\Xa') = p_! L_* (X,\Xa), \]
where $p_!: H_* (X;\rat) \to H_* (X';\rat)$ is the transfer induced by $p$.
\end{thm}
\begin{proof}
We use the notation of the proof of Theorem \ref{thm.signmultfincov}.
Thus $K$ is a simplicial complex triangulating $X$ such that
the Whitney strata of $\Xa$ are triangulated by subcomplexes,
and there is a closed, smooth, oriented $n$-manifold $M$ and a continuous
map $f: M\to X$ such that for a smooth triangulation $L$ of $M$,
\[ [(M,L, \real_M [n], f)] = 2^s [(X,K,\shs, \id_X)] \in S_n^{PL} (X).  \]
For the covering spaces we have
\[ 2^s [(X', K', p^* \shs, \id_{X'})] = 
  [(M', L', \real_{M'} [n], f')] \in S_n^{PL} (X'). \]
By Lemma \ref{lem.flbordinvonsignhom},
\[ L_* (X,\Xa_K) = (\id_X)_* L_* (X,\Xa_K) = 2^{-s} f_* L_* (M) \]
and
\[ L_* (X',\Xa_{K'}) = (\id_{X'})_* L_* (X',\Xa_{K'}) = 2^{-s} f'_* L_* (M'). \]
Using Lemma \ref{lem.vrrsmooth}, the transfer of the manifold
characteristic class is given by
\[  p_{M!} L_* (M) = L_* (M'). \]
The cartesian diagram (\ref{equ.cartdiammprxxpr}) entails
the base change relation $f'_* p_{M!} = p_! f_*$ of type (\ref{equ.basech}).
Therefore,
\begin{align*}
p_! L_* (X,\Xa) 
&= p_! L_* (X,\Xa_K) = 2^{-s} p_! f_* L_* (M) = 2^{-s} f'_* p_{M!} L_* (M) \\
&= 2^{-s} f'_* L_* (M') = L_* (X',\Xa_{K'}) = L_* (X',\Xa').
\end{align*}
\end{proof}
Theorem \ref{thm.lclmult}
on the multiplicativity of L-classes under finite coverings can now
alternatively be deduced as a corollary to Theorem \ref{thm.vrr}
by noting that
\[ p_* L_* (X',\Xa') = p_* p_! L_* (X,\Xa) = d\cdot L_* (X,\Xa). \]

\begin{remark}
Given a topological covering map $q:Y' \to Y$ of finite degree, the methods of this section yield a transfer
\[ q_!: S^{PL}_* (Y) \longrightarrow S^{PL}_* (Y') \]
on signature homology:
For an admissible PL-representative $(X^n, K, \shs, f:X\to Y)$, we
form the pullback cartesian diagram
\[ \xymatrix{
X' = X\times_Y Y' \ar[r]^>>>>>{f'} \ar[d]_p & Y' \ar[d]^q \\
X \ar[r]_f & Y.
} \]
Then $p$ is a topological covering map and
by Lemma \ref{lem.triangofcoverings} there is a triangulation of $X'$ by
a simplicial complex $K'$ such that $p$ is simplicial with respect to $K',K$.
As $p$ is a local homeomorphism, $X'$ is an oriented PL pseudomanifold
which is compact, since $X$ is compact and $p$ has finite degree.
Pullback yields a self-dual sheaf $p^* \shs \in SD(X', \Xa_{K'})$ and we define
the signature homology transfer to be
\[ q_! [(X, K, \shs, f:X\to Y)] = [(X', K', p^* \shs, f':X'\to Y')] \in S^{PL}_n (Y'). \]
Then $q_!$ is well-defined, as an admissible PL-bordism
$(W^{n+1}, Q, \sht, F:W\to Y)$ can similarly be lifted to an admissible
PL-bordism $(W' = W\times_Y Y', Q', p^*_0 \sht, F':W'\to Y')$.
\end{remark}

\section{Varieties with Canonical Singularities}
\label{sec.canonicalsing}

Algebraic varieties are always assumed to be defined over the field of complex numbers.
The \emph{irregularity} $q(X)$ of a normal projective variety $X$ is defined to be
\[ q(X) = \dim H^1 (X;\OO_X). \]
Let $\omega_X$ denote the canonical sheaf and $K_X$ the corresponding
canonical Weil divisor class, $\omega_X = \OO_X (K_X)$. 
One says that $X$ has \emph{canonical singularities}, if the following two conditions
are satisfied:\\
(i) for some integer $r\geq 1,$ the Weil divisor $rK_X$ is Cartier, and \\
(ii) if $f:\widetilde{X} \to X$ is a resolution of $X$ and $\{ E_i \}$ the family
of all exceptional prime divisors of $f$, then
 \[ rK_{\widetilde{X}} = f^* (rK_X) + \sum a_i E_i \text{ with all } a_i \geq 0. \]
According to \cite[Thm. 8.2]{kawamata}, for a normal projective variety $X$
with only canonical singularities, $K_X$ numerically
equivalent to zero ($K_X \equiv 0$) implies that there is a positive integer
$m$ such that $mK_X \sim 0$ (linear equivalence), that is, the double dual of the $m$-th tensor power
of $\omega_X$ is isomorphic to $\OO_X$.

The arguments of the next section rely heavily on a structural analysis
using Albanese maps. Let us briefly review this material with a particular
emphasis on singularities.
First, let $X$ be a smooth connected projective variety. The image of the map
$H_1 (X;\intg) \to H^0 (X;\Omega^1_X)^*$ given by
$z \mapsto (\omega \mapsto \int_z \omega)$, $z\in H_1 (X;\intg),$
$\omega \in H^0 (X;\Omega^1_X),$ is a lattice in $H^0 (X;\Omega^1_X)^*$.
The quotient 
\[ \Alb (X) = \frac{H^0 (X;\Omega^1_X)^*}{H_1 (X;\intg)} \]
is an Abelian variety called the \emph{Albanese variety} of $X$.
Using Serre duality, its dimension is
\[ \dim \Alb (X) = \dim H^0 (X;\Omega^1_X) = \dim H^1 (X; \OO_X) = q(X), \]
the irregularity of $X$.
The Albanese variety comes with a morphism
\[ \alpha_X: X\longrightarrow \Alb (X),~ \alpha_X (x)(\omega) =  \int_{x_0}^x \omega, \]
called the \emph{Albanese map},
where $x_0 \in X$ is a basepoint and one integrates along a path from $x_0$ to $x\in X$.
Modulo the lattice, this is independent of the choice of path.
The Albanese map induces an isomorphism
\[ \alpha_{X*}: H_1 (X;\intg)/\operatorname{Torsion} \stackrel{\cong}{\longrightarrow}
  H_1 (\Alb (X);\intg). \]
If $\kappa (X)=0$, where $\kappa (X)$ denotes the Kodaira dimension of a variety,
then $\alpha_X$ is an algebraic fiber space, that is, it is a surjective morphism with
connected fibers (\cite[Thm. 1]{kawamatacharabvar}).

Given any, possibly singular, connected projective variety $X$, 
let $f:\widetilde{X} \to X$ be any projective
resolution of singularities. The Albanese variety of $X$ is then defined to be
$\Alb (X) := \Alb (\widetilde{X})$. This is independent of the choice of
resolution $\widetilde{X}$, since any two resolutions of $X$ are dominated
by a third one and the Albanese variety as an abelian variety contains
no rational curves. (See e.g. \cite{humengzhang}.) One may define an ``Albanese map''
$X \dashrightarrow \Alb (X)$ as the composition
\[ X \stackrel{f^{-1}}{\dashrightarrow} \widetilde{X}
  \stackrel{\alpha_{\widetilde{X}}}{\longrightarrow} \Alb (\widetilde{X}), \]
but this is in general only a rational map. If one wishes
the property $\dim \Alb (X)=q(X)$ to hold even for singular $X$, then one needs
$q(\widetilde{X}) = q(X)$,
which is not true for arbitrary singularities, but does hold for
$X$ with only \emph{rational} singularities, i.e. if
$R^i f_* \OO_{\widetilde{X}} =0$ for every resolution of singularities
$f:\widetilde{X} \to X$ and for every $i>0$. (It is sufficient to check this for
one resolution of singularities.) This suggests that varieties with only rational singularities
may have a good Albanese map.
Indeed, according to \cite[Lemma 8.1, p. 41]{kawamata},
see also \cite[Prop. 2.3]{reidprojmorkaw},
for a complete normal connected algebraic variety with only rational singularities,
the Albanese map $\alpha_X:X\to \Alb (X)$ is a morphism.
Since canonical singularities are rational, all of the above remarks apply to
varieties with canonical singularities.
Thus if $X$ is a connected projective variety with canonical singularities, then the Albanese map
$\alpha_X: X\to \Alb (X)$ is well-defined, $\dim \Alb (X)=q(X)$, and
if $\widetilde{X} \to X$ is any resolution of singularities, then the diagram
\begin{equation} \label{equ.albxresdiagram}
\xymatrix{
\widetilde{X} \ar[r]^>>>>{\alpha_{\widetilde{X}}} \ar[d] & \Alb (\widetilde{X}) \ar@{=}[d] \\
X \ar[r]^>>>>{\alpha_X} & \Alb (X)
} \end{equation}
commutes. 
Note that if $\kappa (X)=0$, then $\kappa (\widetilde{X})=0$ by
birational invariance of the Kodaira dimension, and thus the above square implies that
$\alpha_X$ is surjective 
and has connected fibers, since $\alpha_{\widetilde{X}}$ is a nonsingular algebraic fiber space.
By \cite[Thm. 8.2]{kawamata},
normal projective varieties with only canonical singularities and $K_X$ numerically
equivalent to zero have $\kappa (X)=0$.

In the following statement of a theorem of Kawamata, we have also incorporated
some remarks of \cite{grebetal}. 
\begin{prop}(Kawamata, \cite[Thm. 8.3]{kawamata}.)   \label{prop.kawamataetalebundle}
Let $X$ be a normal connected projective variety with at worst canonical singularities and 
$K_X$ numerically equivalent to zero. Then the Albanese map
$\alpha_X: X\to \Alb (X)$ is surjective, has connected fibers, and
there exists an abelian variety $E$ and a 
finite degree \'etale cover $E\to \Alb (X)$ such that the fiber product $X \times_{\Alb (X)} E$,
\begin{equation} \label{dia.xalbxe}
\xymatrix{
X \times_{\Alb (X)} E \ar[r] \ar[d] & E \ar[d] \\
X \ar[r]_{\alpha_X} & \Alb (X),
} 
\end{equation}
is isomorphic to a product $X \times_{\Alb (X)} E \cong F\times E$,
where $F$ is a normal projective variety with canonical singularities.
If the canonical divisor of $X$ is trivial, then $F$ has also trivial canonical divisor.
This isomorphism is an
isomorphism \emph{over} $E$, that is, the diagram
\[ \xymatrix{
F\times E \ar[rr]^{\beta}_{\cong} \ar[rd]_{\pi_2} & & X\times_{\Alb (X)} E \ar[ld] \\ 
& E &
} \]
commutes.
\end{prop}
The composition
\[ p:F\times E \cong X\times_{\Alb (X)} E \to X \]
is a degree $d$ \'etale covering, where $d$ is the degree of $E\to \Alb (X)$.
Proposition \ref{prop.kawamataetalebundle}
implies that
\begin{equation} \label{inequ.qxleqdimx}
q(X) = \dim \Alb (X) = \dim E \leq \dim (E\times F) =\dim X 
\end{equation}
for normal projective varieties with only canonical singularities and numerically
trivial canonical divisor. Using the isomorphism $\beta$, we rewrite diagram
(\ref{dia.xalbxe}) as
\begin{equation} \label{dia.kawamcoveringdia}
\xymatrix{
F\times E \ar[r]^{\pi_2} \ar[d]_p & E \ar[d]^{p_E} \\
X \ar[r]_{\alpha_X} & \Alb (X).
} 
\end{equation}
Restricting $p$ to a connected component of a preimage of a connected fiber of $\alpha_X$
yields a homeomorphism between $F$ and that fiber. In particular, $F$ is connected.
Since $\Alb (X)$ is a complex torus, its fundamental group is free abelian.
Therefore, the \'etale cover $E\to \Alb (X)$ is a Galois cover.
Thus the pullback cover $p: F\times E \to X$ is a Galois cover as well.

\section{The Hodge L-Class}

Let $X$ be a complex algebraic variety.
The relative Grothendieck ring $K_0 (\Var/X)$ of complex algebraic
varieties over $X$ has been introduced by Looijenga in \cite{looijenga} 
in the context of motivic integration and
studied by Bittner \cite{bittner}. By variety one means here a separated scheme
of finite type over $\Spec (\cplx)$. As a group, $K_0 (\Var/X)$ is the quotient
of the free abelian group of isomorphism classes of algebraic morphisms
$Y\to X$, modulo the additivity relation generated by
\[ [Y\to X] = [Z\hookrightarrow Y\to X] + [Y-Z\hookrightarrow Y\to X] \]
for $Z\subset Y$ a closed algebraic subvariety of $Y$.
A morphism $f: X' \to X$ induces a  pushforward
\[ f_*: K_0 (\Var/X') \longrightarrow K_0 (\Var/X) \]
by composition:
\[ f_* [Y\to X'] = [Y \longrightarrow X' \stackrel{f}{\longrightarrow} X]. \]
The motivic Hirzebruch natural transformation
\[ T_{y*}: K_0 (\Var /X) \longrightarrow H^{BM}_{2*} (X)\otimes \rat [y], \]
where $H^{BM}_* (-)$ denotes Borel-Moore homology,
has been introduced by Brasselet-Sch\"urmann-Yokura in \cite{bsy}.
The Hirzebruch class of a, generally singular, variety $X$ is defined to be
\[ T_{y*} (X) = T_{y*} ([\id_X]). \]
Borel-Moore homology is covariant with respect to proper maps.
If $f:Y\to X$ is proper, then
\[ T_{y*} [Y \stackrel{f}{\longrightarrow} X] =
 T_{y*} [Y \stackrel{\id}{\longrightarrow} Y \stackrel{f}{\longrightarrow} X] =
 T_{y*} (f_* [\id_Y]) =
 f_* T_{y*} ([\id_Y]) = f_* T_{y*} (Y). \]
For $y=1$, the class $T_{1*} (X)$ is called the 
\emph{Hodge L-class} of $X$, as it is for compact $X$ closely related to, 
but generally different from,
the Goresky-MacPherson L-class $L_* (X)$. 
If $X$ is smooth and pure-dimensional, then
\[ T_{y*} (X) = T^*_y (TX)\cap [X], \]
where $T^*_y (TX)$ is Hirzebruch's cohomological generalized Todd class of the
tangent bundle $TX$ of $X$. For $y=1$, $T^*_1 (TX)=L^* (X),$ the Hirzebruch L-class
of the tangent bundle $TX$. Therefore, for $X$ smooth, compact and pure-dimensional,
\[ T_{1*} (X) = T^*_1 (TX)\cap [X] = L^* (TX)\cap [X] = L_* (X), \]
the Goresky-MacPherson L-class of $X$.

\begin{lemma} \label{lem.tyexc}
Let $E=f^{-1} (x)$ be the exceptional curve of a minimal resolution $f$ of a rational double
point $x$ in a normal surface. If $n$ is the number of irreducible components of $E$, then
\[ T_{y*} (E) = \sum_{i=1}^n T_{y*} (\pr^1_i \to E) - (n-1) T_{y*} (\pt \to E). \]
In particular, for $y=1,$ the degree $0$ component is
\[ T_{1,0} (E) = -(n-1)[\pt]. \] 
\end{lemma}
\begin{proof}
Let $e\in E$ be any point and $j: \{ e \} \hookrightarrow E$ the corresponding
inclusion, a proper map. Then $T_{y*} [j] = j_* T_{y*} (\{ e \})$ and,
since $E$ is path connected, and thus $j_* = j'_*$ for any other inclusion
$j': \{ e' \} \hookrightarrow E,$ the element
$p = T_{y*} [\{ e \} \hookrightarrow E]$ is independent of the choice of $e$.
There are three cases: $x$ is either an $A_n$ ($n\geq 1$), $D_n$ ($n\geq 4$) or an
$E_n$ ($n=6,7,8$) singularity. In any case, $E$ is a tree
$E=\pr^1_1 \cup \cdots \cup \pr^1_n$ of $n$ projective lines.
We establish the $A_n$ case by an induction on $n$.
For $n=1$, we have $E=\pr^1$ and the formula holds. Suppose the formula holds for
$A_{n-1}$ singularities ($n\geq 2$) and $E$ is associated to an $A_n$ singularity. Then
\begin{align*}
T_{y*} (E)
&= T_{y*} (\pr^1_n \cup_{\pt} A_{n-1}) \\
&= T_{y*} [\aff^1_n \to E] + T_{y*} [A_{n-1} \to E] \\
&= T_{y*} [\pr^1_n \to E] - p + \sum_{i=1}^{n-1}T_{y*} [\pr^1_i \to E] - (n-2)p \\
&= \sum_{i=1}^n T_{y*} [\pr^1_i \to E] - (n-1) p.
\end{align*}
For a $D_4$ singularity, $E=\pr^1_4 \cup_{\pt} A_3$ with the intersection point in $\pr^1_2$ of
$A_3 = \pr^1_1 \cup \pr^1_2 \cup \pr^1_3$ and different from the two intersection points in $A_3$.
Thus we can argue as above to get $T_{y*} (E) = \sum_{i=1}^4 T_{y*} [\pr^1_i \to E] - 3p$.
Also the induction step from $D_{n-1}$ to $D_n$ can be carried out as above.
Similarly for $E_n$ singularities.
Now for $y=1$ we obtain, using that $\pr^1$ is smooth,
\[ T_{1*} [\pr^1_i \to E] = \incl_* T_{1*} (\pr^1) = \incl_* L_* (\pr^1) = \incl_* [\pr^1]. \]
This class is concentrated in degree $2$ and has no degree $0$ component. Consequently,
\begin{align*}
T_{1,0} (E) 
&= -(n-1)T_{1,0} [\pt \to E] = -(n-1)\incl_* T_{1,0} (\pt) \\
&= -(n-1)\incl_* L_0 (\pt) = -(n-1)[\pt]. 
\end{align*}
\end{proof}

\section{Proof of the BSY-Conjecture for Trivial Canonical Class}
\label{sec.proofofbsy}

This section proves Theorem \ref{thm.bsy} of the introduction.
Before we can manage $3$-folds, we need to analyze the surface case.
Let $F$ be a normal connected projective surface with at worst canonical singularities
and $K_F$ trivial. By (\ref{inequ.qxleqdimx}), such a
surface has irregularity $q(F) \in \{ 0,1,2 \}$. The most interesting situation
is $q(F)=0$, which we shall investigate first.\\

Let $F$ be a normal connected projective surface with at worst canonical singularities and
trivial canonical class.
Since $F$ is normal, it has no singularities in complex codimension $1$.
Thus $F$ has only isolated singularities. As $F$ is compact, there are only finitely
many of them, call them $x_1,\ldots, x_m\in F$.
Two-dimensional canonical singularities
are precisely the du Val singularities.
In dimension $2,$ crepant resolutions of du Val singularities always exist and are unique.
(Moreover, every crepant resolution of a surface is minimal.)
Let $f: F_0 \to F$ be the unique crepant resolution of $F$.
Then 
\[ K_{F_0} = f^* K_F =0. \]
Since by definition, $\Alb (F_0)=\Alb (F)$,
we have
\[ q(F_0) = \dim \Alb (F_0) = \dim \Alb (F) = q(F). \]
We summarize: The surface $F_0$ is smooth, connected, has trivial canonical class and $q(F_0)=q(F)$.
So if $q(F)=0$, then $F_0$ is a smooth $K3$ surface. All such surfaces are diffeomorphic to
each other and have Hodge numbers
\[ h^{2,0} (F_0)=h^{0,2} (F_0)=1,~ h^{1,1} (F_0)=20. \]
By the Hodge index theorem, the signature $\sigma (F_0)$ of $F_0$ is
\[ \sigma (F_0)= 2h^{2,0} - h^{1,1} +2=4-20=-16. \]

We return to $F$ with unrestricted irregularity.
Let $E_i = f^{-1} (x_i)$ be the exceptional curve over the singular point $x_i$ and let
\[ E_i = E_{i1} \cup \cdots \cup E_{id_i} \]
be its decomposition into irreducible components.
Then the $E_{ij}$ are rational $(-2)$-curves, i.e. $E_{ij}$ is a $\pr^1$ with
self-intersection $E^2_{ij} =-2$.
Two curves $E_{ij}$ and $E_{ik}$ intersect transversely in a single point if they are not disjoint,
and no three of them intersect. 
The dual graph of the exceptional divisor $E_i$ is a simply laced
Dynkin diagram
$A_n$ ($n\geq 1$), $D_n$ ($n\geq 4$), or $E_n$ ($n=6, 7, 8$), where $n=d_i$
is the number of irreducible components of $E_i$.
Let $T_i \subset F_0$ be a small regular neighborhood of $E_i,$ $i=1,\ldots, m,$
such that $T_i \cap T_j =\varnothing$ for $i\not= j$.
Then $T_i$ is a compact manifold with boundary which contains $E_i$ as a 
deformation retract. In fact, $T_i$ is diffeomorphic to the (abstract) manifold $P_i$
obtained by plumbing on the normal disc bundles of the $E_{ij}$ according to the
dual graph (the Dynkin diagram). Thus the homology of $P_i$ is the homology of $E_i$;
in particular
\[ H_2 (T_i;\real) \cong H_2 (P_i;\real)\cong H_2 (E_i;\real) \cong \real^{d_i}, \]
generated by the fundamental classes $[E_{ij}]$.
Moreover, the intersection form on $P_i$ 
(and hence the one on $T_i$)
is given by the intersection matrix of the dual
graph of $E_i$, which is the negative of the Cartan matrix of the relevant Dynkin diagram.
The boundary $\partial T_i \cong \partial P_i$ is the link $L_i \subset F$ of the
singular point $x_i$. We can assume that the $T_i$ have been arranged in such a 
way that for suitable closed cone neighborhoods $N_i$ of $x_i$ in $F$, $N_i \cong \cone (L_i),$
$f^{-1} (N_i)=T_i$ and the restriction of $f$ to the complement $M_0$ of the interiors of the
tubes is an isomorphism onto the complement $M$ of the interiors of the $N_i$.

\begin{lemma} \label{lem.surfacerathommfd}
The surface $F$ is a rational homology manifold.
\end{lemma}
\begin{proof}
The statement is equivalent to $H^1 (f^{-1} (x_i);\rat)=0$ for all $i$.
But the latter follows from $H^1 (E_{ij};\rat) = H^1 (\pr^1;\rat)=0$ together
with the fact that every exceptional curve $E_i$ is a tree of $\pr^1$s.
\end{proof}

If $(M,\partial M)$ is a compact oriented manifold with boundary, 
let $\sigma (M,\partial M)$ denote its Novikov signature, that is,
the difference between the number of positive and negative eigenvalues of the intersection
form on the middle homology of $M$, which need not be non-degenerate.
By \cite{siegel}, $\sigma (M,\partial M)$ agrees with the Goresky-MacPherson
signature of the closed pseudomanifold obtained from $(M,\partial M)$ by
coning off the boundary $\partial M$.
By a theorem of Mumford and Grauert (see \cite[Theorem III.2.1]{barthetal})
the intersection matrix of $E_i$ is negative definite. It follows that
\[ \sigma (T_i, \partial T_i)= -d_i. \]

Let $M_0$ be the $4$-manifold (with boundary) obtained by taking the complement of the interior of the 
tubes $T_1 \cup \cdots \cup T_m$.
By Novikov additivity of the signature,
\[ \sigma (M_0,\partial M_0) + \sum_{i=1}^m \sigma (T_i,\partial T_i) = \sigma (F_0). \]
Let $M$ be the $4$-manifold obtained by taking the complement in $F$ of the interior of  
$N_1 \cup \cdots \cup N_m$.
Novikov additivity of the Goresky-MacPherson signature for general Witt spaces
has been established in \cite{siegel} and implies
\[ \sigma (M,\partial M) + \sum_{i=1}^m \sigma (N_i,\partial N_i) = \sigma (F). \]
Since $f$ restricts to an isomorphism $(M_0,\partial M_0)\cong (M,\partial M)$,
we have $\sigma (M_0,\partial M_0) = \sigma (M,\partial M)$.
The signature $\sigma (N_i, \partial N_i)$ equals the signature of $\cone (L_i)$ with
the boundary $L_i$ coned off, that is, the signature of the suspension of the link.
But the signature of any suspension is zero, since the suspension has an orientation-reversing
automorphism given by flipping the suspension coordinate. We arrive at the formula
\begin{equation} \label{equ.signfgen} 
\sigma (F) = \sigma (F_0) + \sum_{i=1}^m d_i. 
\end{equation}
If $q(F)=0$, then this specializes to:
\begin{equation} \label{equ.signf} 
\sigma (F) = -16 + \sum_{i=1}^m d_i. 
\end{equation}
According to \cite[Thm. 18.21]{fletcher},
\[  \sum_{i=1}^m d_i \leq 19. \]
Using (\ref{equ.signf}), we obtain the bound
\[  -16 \leq \sigma (F) \leq 3 \]
for the signature of $F$ with $q(F)=0$.

\begin{lemma}  \label{lem.duvalsurft1fislf}
Let $F$ be a normal connected projective surface with at worst du Val singularities,
trivial canonical divisor and $q(F)=0$. Then
\[ \sigma (F) \in \{ -16, -15,\ldots, 2,3 \} \]
and the Hodge L-class of $F$ equals the Goresky-MacPherson
L-class of $F$:
\[ T_{1*} (F) = L_* (F). \]
\end{lemma}
\begin{proof}
The signature statement has already been established above.
Let $S = \{ x_1,\ldots, x_m \} \subset F$ be the singular set of $F$ and
let $\Ex = \bigcup_{i=1}^m E_i \subset F_0$ be the exceptional set of the unique crepant
resolution $f:F_0 \to F$.
We shall use the commutative diagrams
\begin{equation} \label{equ.f0mexfmsinf0f}
\xymatrix{
F_0 - \Ex \ar@{^{(}->}[r] \ar[d]_{f|}^\cong & F_0 \ar[d]^f \\
F - S \ar@{^{(}->}[r] & F,
} 
\end{equation}
and
\begin{equation} \label{equ.exsinf0f}
\xymatrix{
\Ex \ar@{^{(}->}[r]^{\incl} \ar[d]_{f|} & F_0 \ar[d]^f \\
S \ar@{^{(}->}[r]^{\incl} & F.
} 
\end{equation}
Using additivity of the Hodge L-class and diagram (\ref{equ.f0mexfmsinf0f}), we obtain for
the Hodge L-class:
\begin{align*}
T_{1*} (F)
&= T_{1*} [F-S \hookrightarrow F] + T_{1*} [S \hookrightarrow F] \\
&= T_{1*} [F_0 -\Ex \stackrel{\cong}{\longrightarrow} F-S \hookrightarrow F] + \incl_* T_{1*} (S) \\
&= T_{1*} [F_0 -\Ex \hookrightarrow F_0 \stackrel{f}{\longrightarrow} F] + m[\pt] \\
&= f_* T_{1*} [F_0 -\Ex \hookrightarrow F_0] +m[\pt].
\end{align*}
By additivity on $F_0$,
\[ T_{1*} (F_0) = T_{1*} [F_0 -\Ex \hookrightarrow F_0] + T_{1*} [\Ex \hookrightarrow F_0]. \]
On the other hand, as $F_0$ is smooth (and pure-dimensional),
\[ T_{1*} (F_0) = L_* (F_0) = \sigma (F_0)[\pt] + [F_0] = -16 [\pt] + [F_0]. \]
Also, as $f$ is a degree $1$ map, we have $f_* [F_0] = [F]$ for the fundamental classes.
Therefore, using diagram (\ref{equ.exsinf0f}),
\begin{align*}
T_{1*} (F)
&= f_* (T_{1*} (F_0) - T_{1*} [\Ex \hookrightarrow F_0] )+m[\pt] \\
&= f_* T_{1*} (F_0) - f_* \incl_* T_{1*} (\Ex) + m[\pt] \\
&= f_* T_{1*} (F_0) - \incl_* f|_* T_{1*} (\Ex) + m[\pt] \\
&= f_* (-16[\pt] + [F_0]) - \incl_* f|_* T_{1,0} (\Ex) + m[\pt] \\
&= (m-16)[\pt] + [F] - \incl_* f|_* T_{1,0} (\Ex).  
\end{align*}
By Lemma \ref{lem.tyexc},
\[ T_{1,0} (E_i) = -(d_i -1)[\pt]. \]
Thus, as $F_0$ is path connected,
\[ T_{1,0} (\Ex) = \sum_{i=1}^m T_{1,0} [E_i \hookrightarrow \Ex]= \sum_{i=1}^m (1-d_i)[\pt_i],~ 
  \pt_i \in E_i. \]
Therefore, using the signature formula (\ref{equ.signf}),
\begin{align*}
T_{1*} (F)
&= (m-16)[\pt] + [F] - \incl_* f|_* \sum_{i=1}^m (1-d_i)[\pt_i] \\
&= (m-16)[\pt] + [F] - \left( m- \sum_{i=1}^m d_i \right)[\pt] \\
&= (-16 +  \sum_{i=1}^m d_i)[\pt] + [F] = \sigma (F)[\pt] + [F] \\
&= L_* (F).
\end{align*}
\end{proof}

We proceed to the case $q(F)>0$.
\begin{lemma}  \label{lem.duvalsurfqf12t1fislf}
Let $F$ be a normal connected projective surface with at worst du Val singularities,
trivial canonical divisor and $q(F)>0$. Then
\[ \sigma (F) =0, \]
$F$ is nonsingular and thus the Hodge L-class of $F$ equals the Goresky-MacPherson
L-class of $F$:
\[ T_{1*} (F) = L_* (F) = [F], \]
where $[F]\in H_4 (F;\rat)$ is the fundamental class.
\end{lemma}
\begin{proof}
The surface $F$ has a 
Kawamata covering diagram of type (\ref{dia.kawamcoveringdia}),
\begin{equation} \label{dia.kawamatf} 
\xymatrix{
F' \times E' \ar[r]^{\pi_2} \ar[d]_{p'} & E' \ar[d] \\
F \ar[r]_{\alpha_F} & \Alb (F),
} 
\end{equation}
where the vertical maps are covering projections of equal finite degree.
There are two possibilities: either $q(F)=1$ or $q(F)=2$.
Suppose that $q(F)=1$.
Since $\dim E' = \dim \Alb (F) = q(F) =1,$ $E'$ is a curve. Thus $F'$ is a curve
as well and $\sigma (F' \times E')=0$.
By multiplicativity of the signature under finite covers (Corollary \ref{cor.signmultalg}),
$\sigma (F)=0$.
Let $f:F_0 \to F$ be the unique crepant resolution, $K_{F_0} =0$.
Then $q(F_0)=q(F)=1$ and we may apply the same argument as above to 
$F_0$ instead of $F$ to conclude that $\sigma (F_0)=0$.
But then formula (\ref{equ.signfgen}) shows that $\sum_{i=1}^m d_i =0$,
that is, $d_i =0$ for all $i$. Thus $f$ is an isomorphism and $F$ is nonsingular.

Suppose that $q(F)=2$.
The fibers of $\alpha_F$ in diagram (\ref{dia.kawamatf}) are connected.
Since these are homeomorphic to $F'$, the space $F'$ is connected as well.
As $\dim E' = \dim \Alb (F)=q(F)=2$, the variety $E'$ is a surface and $F'$
is zero-dimensional. Since $F'$ is connected, it is a point and thus 
$F' \times E' = E'$. The map $\pi_2$ is the identity.
The Albanese variety $\Alb (F)$ is a complex torus homeomorphic to
$S^1 \times S^1 \times S^1 \times S^1$ and thus has signature $\sigma (\Alb (F))=0$.
By Corollary \ref{cor.signmultalg}, $\sigma (E')=0$, and, applying that corollary
to the finite covering $p': E' \to F,$ the signature of $F$ vanishes.
By the same token, $\sigma (F_0)=0$ for the unique crepant resolution
$f: F_0 \to F$ of $F$. It then follows as in the case $q(F)=1$
that $F$ is nonsingular.
\end{proof}

Let $X$ be a normal connected projective $3$-fold 
with only canonical singularities, with trivial canonical class and $q(X)>0$.
Then $X$ has a 
Kawamata covering diagram (\ref{dia.kawamcoveringdia}),
\[ \xymatrix{
F \times E \ar[r]^{\pi_2} \ar[d]_p & E \ar[d]^{p_E} \\
X \ar[r]_{\alpha_X} & \Alb (X),
} \]
where $p$ and $p_E$ are covering projections of finite degree $d$.
The variety $F$ is connected, being homeomorphic to the fiber of $\alpha_X$.
By (\ref{inequ.qxleqdimx}), $q(X) \leq \dim X.$
Thus, as $X$ has dimension $3$, there are the three cases $q(X)\in \{ 1,2,3 \}$.

We begin with the most interesting case $q(X)=1$.
Then $\dim \Alb (X)=q(X)=1$ and both $\Alb (X)$ and $E$ are curves.
Consequently, $F$ is a surface (which is normal, projective, has at worst
canonical singularities and trivial canonical class).
\begin{remark}
The above covering diagram together with Lemma \ref{lem.surfacerathommfd}
implies that $X$ is a rational homology manifold.
\end{remark}
Using cartesian multiplicativity of the Goresky-MacPherson L-class (Proposition 
\ref{prop.lclasscartesianmult}),
\[ L_* (F\times E) = L_* (F)\times L_* (E) =
  (\sigma (F)[\pt_F] + [F])\times [E] = \sigma (F)\cdot [\pt_F \times E] + [F\times E]. \]
Thus
\begin{align*} 
p_* L_* (F\times E) 
  &= p_* ( \sigma (F)\cdot [\pt_F \times E] + [F\times E]) \\
 &= \sigma (F)\cdot p_* [\pt_F \times E] + d [X].
\end{align*}
On the other hand by Corollary \ref{cor.lclassmultalgcover},
\[ p_* L_* (F\times E) = d L_* (X). \]
Therefore,
\[ L_* (X) = \frac{\sigma (F)}{d} \cdot p_* [\pt_F \times E] + [X]. \]
We shall now carry out a parallel computation for the 
motivic Hodge L-class $T_{1*} (X)$.
An exterior product $\times$ on relative Grothendieck groups of varieties is given by
\[ 
\begin{array}{rcl}
K_0 (\Var/X) \times K_0 (\Var/Y) & \longrightarrow & K_0 (\Var/ X\times Y) \\
\mbox{[} f:Z\to X \mbox{]} \times \mbox{[} g:Z'\to Y \mbox{]} 
   & \mapsto & \mbox{[} f\times g: Z\times Z' \to X\times Y \mbox{]}.
\end{array}
\]
According to \cite{schueryokexpo}, the class $T_{y*}$ commutes with this
exterior product, that is, the diagram
\[ \xymatrix{
K_0 (\Var /X) \times K_0 (\Var /Y) \ar[r]^{\times} \ar[d]_{T_{y*} \times T_{y*}} &
   K_0 (\Var /X\times Y) \ar[d]^{T_{y*}} \\
H^{BM}_{2*} (X) \otimes \rat[y] \times H^{BM}_{2*} (Y) \otimes \rat[y] \ar[r]^>>>>>{\times}  &
   H^{BM}_{2*} (X\times Y)  \otimes \rat[y]
} \]
commutes. Thus
\[ T_{1*} (F\times E) = T_{1*} (F)\times T_{1*} (E). \]
Since $E$ is smooth, we have $T_{1*} (E)=L_* (E)$, and by 
Lemmata \ref{lem.duvalsurft1fislf}, \ref{lem.duvalsurfqf12t1fislf},
$T_{1*} (F) = L_* (F)$. Therefore,
\[ T_{1*} (F\times E) = L_* (F)\times L_* (E) = L_* (F\times E) \]
and consequently
\[  p_* T_{1*} (F\times E) = p_* L_* (F\times E) = \sigma (F)\cdot p_* [\pt_F \times E] +d[X]. \]
We shall now require an analog of 
Corollary \ref{cor.lclassmultalgcover} for $T_{1*}$. This can readily be
deduced from a more general result of Cappell, Maxim, Schürmann and Shaneson,
\cite[Theorem 5.1]{cmssequivcharsing}, which asserts that if a finite group $G$
acts algebraically on a complex quasi-projective variety $X'$, then
\begin{equation} \label{equ.cmss}
T_{y*} (X'/G) = \frac{1}{|G|} \sum_{g\in G} \pi^g_* T_{y*} (X';g), 
\end{equation}
where $\pi^g: X'^g \to X'/G$ is the composition of the inclusion $X'^g \hookrightarrow X'$
of the fixed point set $X'^g$ and the projection $\pi:X' \to X'/G$. The classes $T_{y*} (X';g)$
are Hodge-theoretic Atiyah-Singer classes in Borel-Moore homology,
supported on the fixed point sets,
\[ T_{y*} (X';g) \in H^{\operatorname{BM}}_{2*} (X'^g) \otimes \cplx [y]. \]
Formula (\ref{equ.cmss}) implies in particular:
\begin{prop} \label{prop.tymultfincov}
If $\pi: X' \to X$ is a finite degree $d$ Galois covering with $X'$ complex quasi-projective,
then
\[ \pi_* T_{y*} (X') = d T_{y*} (X). \] 
\end{prop}
One can also obtain that relation from the Verdier Riemann-Roch formula of
\cite[Cor. 3.1.3]{bsy}, noting that the vertical tangent bundle $T_p$ of a covering
map $p$ is the rank zero bundle.
Now, as $F\times E$ is projective and the \'etale covering 
$p: F\times E \to X$ is indeed a Galois covering, Proposition \ref{prop.tymultfincov}
implies that 
$p_* T_{1*} (F\times E) = d T_{1*} (X),$ and thus
\[ T_{1*} (X) = \frac{1}{d} p_* T_{1*} (F\times E) =  \frac{1}{d} p_* L_{*} (F\times E) = L_* (X). \]
This proves the Brasselet-Sch\"urmann-Yokura conjecture for 
the case $q(X)=1$.\\

Suppose that $q(X)=2$.
Then $\Alb (X)$ and $E$ are surfaces and $F$ is a (normal) curve.
As $E$ is an abelian variety, it
is a complex torus $\cplx^2/L$ and thus $\sigma (E)=0$.
It follows that 
\[ L_* (F\times E) = L_* (F)\times L_* (E) =
  [F] \times (\sigma (E)[\pt_E] + [E]) 
   = [F\times E]. \]
Since $F$ is a normal curve, it is nonsingular and therefore
$T_{1*} (F) = L_* (F) = [F]$. Since $E$ is smooth as well,
\[ T_{1*}(F\times E) = T_{1*} (F)\times T_{1*} (E) = L_* (F)\times L_* (E) = L_* (F\times E). \]
Consequently, 
\[ T_{1*} (X) = L_* (X) = \frac{1}{d} p_* [F\times E] = [X]. \]

Finally, suppose that $q(X)=3$.
The fibers of $\alpha_X$ are connected.
Since these are homeomorphic to $F$, the space $F$ is connected as well.
As $\dim E = \dim \Alb (X)=q(X)=3$, the variety $E$ is a $3$-fold and $F$
is zero-dimensional. Since $F$ is connected, it is a point and thus 
$F \times E = E$. The map $\pi_2$ is the identity.
The variety $E$ is a complex torus, so in particular nonsingular.
Therefore, $T_{1*} (E) = L_* (E) = L^* (T^6)\cap [E] = 1\cap [E]= [E].$ Applying $p_*$ and dividing by $d$, we arrive
at $T_{1*} (X) = L_* (X) = [X]$.
The BSY-conjecture is thus proved for all normal, connected, projective $3$-folds $X$
with only canonical singularities, trivial canonical class, and $q(X)>0$.\\

Obviously, our methods give some information in dimensions higher than $3$
as well. For example, if $X$ is a normal connected projective $4$-fold with canonical singularities,
trivial canonical class and $q(X)=1$, then $T_{1*}(X)=L_* (X)$ provided
$q(F)>0$ in the Albanese-Kawamata covering diagram of $X$.
A general inductive proof for arbitrary dimension cannot be put together from
the methods presented here, since the condition $q(X)>0$ does not propagate
inductively.

\section{Examples and Realization}

We discuss briefly some examples of projective $3$ folds with 
canonical singularities, trivial canonical class and $q(X)>0$.
In the previous section we have seen that if $q(X)>1$, then such an 
$X$ is in fact a topological manifold. Therefore, we shall focus on $q(X)=1$.

A projective surface $F$ is called a (possibly singular)
\emph{$K3$ surface} if it has only du Val singularities,
$\omega_F \cong \OO_F$ and $q(F)=0$ (\cite[Section 2.1]{brown}).
Let $a_0, \ldots, a_n$ be positive integers and let
$\pr (a_0,\ldots, a_n)$ denote the associated weighted projective space.
We denote by
$F_{d_1,\ldots, d_c} \subset \pr (a_0,\ldots, a_n)$ a general element of the family of all weighted
complete intersections of multidegree $\{ d_1, \ldots, d_c \}$.
Such a weighted complete intersection $F_{d_1,\ldots, d_c}$ has dimension $n-c$.
All weighted complete intersections (and weighted projective spaces) are $V$-manifolds,
i.e. locally quotients of $\aff^n$ by a finite group action.
Several lists of $K3$ surfaces embedded in weighted projective space exist.
The first of these was Reid's list of ``famous'' $95$ weighted hypersurfaces of \cite{reidfamous95}.
Independently, Yonemura \cite{yonemura} gave a list $W_4$ of $95$ weights which
coincides with Reid's set $A'_4$.
Then $84$ families of $K3$ surfaces in codimension $2$ were listed by
Fletcher \cite[Section 13.8]{fletcher}. Alt{\i}nok found $70$ families in codimension $3$
and $142$ families in codimension $4$, \cite{altinokthesis}.
There exists a computer-supported database of candidate polarised $K3$ surfaces compiled
by Gavin Brown, see \cite{brown}.
The integers
\[ \{ -16, -15, -14, -13, -11, -10, -9, -8, -7, -6, -5, -4, -3, -2, -1, 0, 1, 2 \} \]
are realized as signatures $\sigma$ of weighted hypersurfaces by members of Reid's 95.
Note that $-12$ is not realized in codimension $1$, but is realized as the signature of
weighted codimension $2$ complete intersections by members of Fletcher's 84;
see the table for examples that realize each integer $-16,\ldots, 2$:
\begin{table}
\[ \begin{array}{l|l|r}
\text{Weighted } K3 & \text{Singularities} & \sigma (F) \\ \hline
F_4 \subset \pr (1, 1, 1, 1) &       & -16 \\ \hline 
F_5 \subset \pr (1, 1, 1, 2) & A_1 & -15 \\ \hline 
F_8 \subset \pr (1, 1, 2, 4) & 2A_1 & -14 \\ \hline 
F_6 \subset \pr (1 , 1, 2, 2) & 3A_1 & -13 \\ \hline 
F_{4,4} \subset \pr ( 1,1 ,2 ,2,2 ) & 4A_1 & -12 \\ \hline 
F_{10} \subset \pr (1, 2,2 ,5 ) & 5A_1 & -11 \\ \hline 
F_{8} \subset \pr (1, 2, 2, 3) & 4A_1, A_2 & -10 \\ \hline 
F_{9} \subset \pr (1, 2,3 ,3 ) & A_1, 3A_2 & -9 \\ \hline 
F_{16} \subset \pr (1, 3, 4, 8) & A_2, 2A_3 & -8 \\ \hline 
F_{12} \subset \pr (1, 3, 4, 4) & 3A_3 & -7 \\ \hline 
F_{12} \subset \pr (2, 2, 3, 5) & 6A_1, A_4 & -6 \\ \hline 
F_{12} \subset \pr (2, 3, 3, 4) & 3A_1, 4A_2 & -5 \\ \hline 
F_{14} \subset \pr (2, 3, 4, 5) & 3A_1, A_2, A_3, A_4 & -4 \\ \hline 
F_{15} \subset \pr (2, 3, 5, 5) & A_1, 3A_4 & -3 \\ \hline 
F_{18} \subset \pr (3, 4, 5, 6) & A_1, 3A_2, A_3, A_4 & -2 \\ \hline 
F_{19} \subset \pr (3, 4, 5, 7) & A_2, A_3, A_4, A_6 & -1 \\ \hline 
F_{24} \subset \pr (3, 4, 7, 10) & A_1, A_6, A_9 & 0 \\ \hline 
F_{25} \subset \pr (4, 5, 7, 9) & A_3, A_6, A_8 & 1 \\ \hline 
F_{30} \subset \pr (5, 6, 8, 11) & A_1, A_7, A_{10} & 2 \\ \hline 
\end{array} \]
Signatures realized by $K3$ surfaces with singularities.
\end{table}
It is apparently not known at the time of writing whether the remaining theoretically possible
integer $3$ is in fact the signature of a $K3$ surface with singularities.

Let $E$ be an elliptic curve and let $F$ be any $K3$ surface as above, for example one
of the weighted complete intersections listed.
Then $E$ and $F$ both have trivial canonical sheaf.
Thus, the $3$-fold $X = F\times E$ has trivial canonical sheaf. 
Let $f:F_0 \to F$ be the unique crepant resolution of $F$.
Then $g=f\times \id_E: F_0 \times E \to F\times E$ is a projective resolution of singularities 
and thus $\Alb (F\times E)=\Alb (F_0 \times E)$.
Since $F_0 \times E$ is smooth,
the Albanese map induces an isomorphism
\[ H_1 (F_0 \times E;\intg)/\operatorname{Torsion} \stackrel{\cong}{\longrightarrow}
  H_1 (\Alb (F_0 \times E);\intg). \]
Now $H_1 (F_0 \times E;\intg) = H_1 (E;\intg) = \intg^2,$
so that 
$H_1 (\Alb (F_0 \times E);\intg) = \intg^2.$
As $\Alb (F_0 \times E)$ is a complex torus $\cplx^k/L$ for some $k$, and
$H_1 (\cplx^k/L;\intg)=\intg^{2k}$,
it follows that $k=1$ and thus 
\[ q(F\times E) = \dim \Alb (F\times E)= \dim \Alb (F_0\times E)=k=1. \]
Finally, $X=F\times E$ has canonical singularities, since the resolution $g$
is crepant; in fact, $K_{F_0 \times E} =0 = g^* K_{F\times E}.$

We remark on the side that smooth complex projective $3$-folds $X$ with trivial canonical bundle, 
$q(X)=0$ but infinite fundamental group
have been completely classified by Kenji Hashimoto and
Atsushi Kanazawa in \cite{hashimotokanazawa}.
It turns out that these are rather scarce:
Up to deformation equivalence, there are precisely $14$ smooth Calabi-Yau $3$-folds
with infinite fundamental group; $6$ of them are covered by an abelian $3$-fold,
and the remaining $8$ are covered by the product of a smooth $K3$ surface $F$ and an
elliptic curve $E$. The simplest example of the latter type is the
Enriques Calabai-Yau $3$-fold
$X= (F\times E)/\langle (\iota,-1_E)\rangle$, where
$\iota$ is an Enriques involution (i.e. a fixed point free involution)
and $-1_E$ is negation. The fundamental group $\pi_1 (X)$ is a semidirect product
$\intg^2 \rtimes C_2$, where the cyclic group $C_2$ of order two is the Galois group
of the cover generated by the involution and the $C_2$-action on $\intg^2$
is the action of the Galois group on $\pi_1 (E) \cong \pi_1 (F\times E)$.
The abelianization $H_1 (X;\intg)$ is the torsion group $(\intg/_2)^3$.
In particular, $H_1 (X;\intg)$ is trivial modulo torsion, which is consistent with
the Albanese variety $\Alb (X)$ being a point.

\end{document}